\documentclass[10pt]{article}
\usepackage[dvips]{graphicx}
\usepackage{amsmath}
\usepackage{array}
\begin{document}

\bf\title{Geometry of the Riemann Zeta Function}
\author{George H. Nickel}
\maketitle

\tableofcontents

\begin{abstract}
\bf

The Riemann zeta function $\zeta (s)$, where $s = \sigma +it$ and  $(\sigma,t) \in \mathcal{R}$, is related to $\sum_1^{\infty}n^{-s}$ by analytic continuation \cite{edwards}.  The traditional analysis pioneered by Riemann uses complex integration of the function $(-x)^s/(e^x-1)dx/x$ around specific paths.  A geometrical examination of the accumulated steps $n^{-s}$ of the sum reproduces the results of the complex analysis and gives new insights into analytic continuation and the occurrence of zeros.  In particular, there is a detailed symmetry under which each step $n < n_p = [\sqrt{t/2\pi}]$ corresponds to a region of conjugate steps $\bar{n} \in  (n_p^2/(n +1/2), n_p^2/(n -1/2))$. End points of steps prior to the center of symmetry are conjugate to steps having an angle change of an odd multiple of $\pi$, leading to conjugate regions of linked (Euler, or Cornu) spirals.  The sum of the steps in each conjugate region (with end corrections) equals $Q(s)n^{s-1}$, where $Q(s)=n_p^{1-2s}e^{i(t+\pi/4)}$, independent of $n$.  Analytically, this corresponds to $(2\pi i)^s/\Pi(s-1)$, the coefficient of the ratio of the $nth$ step to the integral of the $nth$ imaginary axis pole of the generating function. The functional equation of $\zeta(s)$ is an expression of the equality of the sums over $n$ and the corresponding conjugate regions.  Two points in the complex plane are specifically recognized; a center of symmetry $P(s)$, which is conjugate to itself, and the point conjugate to the origin.  The latter is shown to agree with the previous methods of analytic continuation for $\zeta(s)$, including Riemann's general description of continuation as ``remaining valid'' where the sum diverges.  The sum of $P(s)$ and its symmetric counterpart $Q(s)P(1-s)$ is the Riemann-Siegel (R-S) equation for $\zeta (s)$, and the lowest order geometrically determined difference between the discrete sum $\sum_1^{n_p} n^{-s}$ and $P(s)$ yields the first-order R-S correction as found by Riemann.  A geometrical construction similar to the \textit{lima\c con} shows that zeros occur when the directions of $P(s)$ and $Q(s)P(1-s)$ are opposed and their magnitudes are equal.  The magnitudes of these surfaces are seen to be equal when $\sigma=1/2$, and their intersection at $\sigma = 1/2$ is generally transverse.  Zeros for other values of $\sigma$ require exceptional conditions, such as simultaneous zeros of the sums for $P(s)$ and $P(1-s)$.

An expression is developed for the angles between adjacent steps, \textit{modulo} $2\pi$. Due to the modulus function, these angles appear to be unrelated for large $t/n$, passing through a decreasing succession of odd and even multiples of $\pi$ as $n$ increases.  The shape, orientation and extent of these characteristic patterns are examined systematically, beginning at the symmetry center and progressing outward to the origin and its conjugate point $\zeta(s)$. Analytical results are summarized for comparisons with these features.  In addition to the results cited above, this visualization provides a new approach to investigation of other features, such as the location, number and spacing of zeros, average values, determination of bounds and anomalous behavior near closely spaced zeros.  Use of the Landau formula relating the primes to the sum of $x^\rho$ over roots $\rho = 1/2 + i\alpha_i$ \cite{landau}, shows that the sum of the prime steps $p^{-s}$ accounts, on average, for the decrease of the real component from its initial value of unity to zero for the zeros of $\zeta(s)$.  The geometric methods developed here may also be applicable to generalized $\zeta$-functions, although the detailed symmetry is no longer present.
\end{abstract}
\rm

\section{Introduction}

\begin{figure*}[p]
\centering
\includegraphics[angle=90,width=1.2\textwidth]{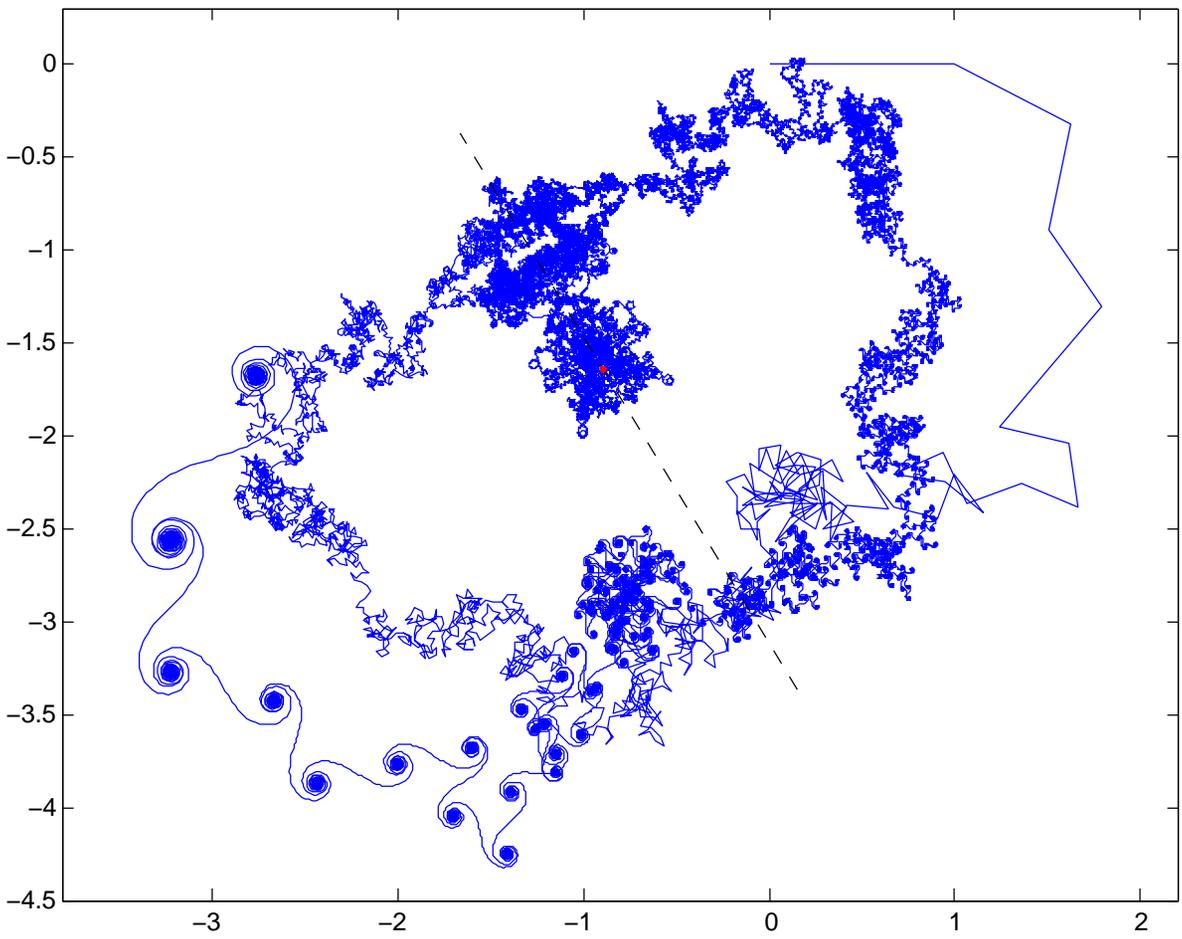}
\caption{Argand Diagram for $\zeta(1/2+ 10^9i)$. }
\label{xy9}
\end{figure*}

\rm
The diagram in Figure \ref{xy9} is a cumulative sum of over 300 million steps $1/n^s$ in the complex plane for $s=1/2+10^9 i$, beginning with $n=1$ at the origin.  Perhaps surprisingly for such a simple prescription, the resulting figure is highly symmetric.  There is a central step at $n=12615$, shown in red. Also indicated is an axis of bilateral symmetry at the angle $121.275^o/301.275^o$.  The symmetry in the scale of opposing features is unique to the value $\sigma=1/2$, and the angles are determined by $t$, which was chosen arbitrarily for this figure.  This symmetry is due to several elementary mathematical effects; the functional properties of the logarithm, the nature of $tlog(n)$ (mod$ ,2\pi$) for large argument, and the addition of steps in the complex plane. Careful inspection reveals that while features on opposite sides of the axis have the same overall dimensions, their detailed structure is not identical.  An extreme example of this is that each of the initial steps is replicated precisely by the center-to-center distances and angles of the large, linked spirals.

 With appropriate attention to convergence issues, the sum of these steps gives one value of the function $\zeta(s)$. In this paper, simple geometric methods involving points, line segments and angles in the complex plane are used to derive the main features of $\zeta(s)$. The results of this detailed yet elementary \textit{geometric number theory} can be correlated with most, if not all, of those of the conventional analysis.  This view extends mathematical understanding and presents new interpretations; in particular, a method of convergence is shown and justified mathematically, in which the truncation of the step sum with a specified end correction provides the properly converged (mathematically termed \textit{continued}) value of $\zeta$. For values of $s$ at which $\zeta=0$, the prime steps are shown to be those for which the value of 1 after the first step eventually decreases to 0. Other examples are explored by this technique, such as the occurrence of zeros, and there are undoubtedly many opportunities for further discoveries.

 Following an introductory discussion of the features of $\zeta (s)$ and its relation to the prime numbers, the remainder of the paper is divided into two parts.  First, calculational methods to treat the geometry of the Argand diagram are presented, and its characteristic patterns are explored systematically.  A selective summary of the results of analytic number theory is then provided to correlate the geometrical results with those of the conventional analytical methods.  This combination of the two approaches adds greatly to the understanding of Riemann's famous function.

\section{Why is $\zeta (s)$ of Interest?}

\begin{figure}[pt]
\centering
\includegraphics[angle=0,width=.6\textwidth]{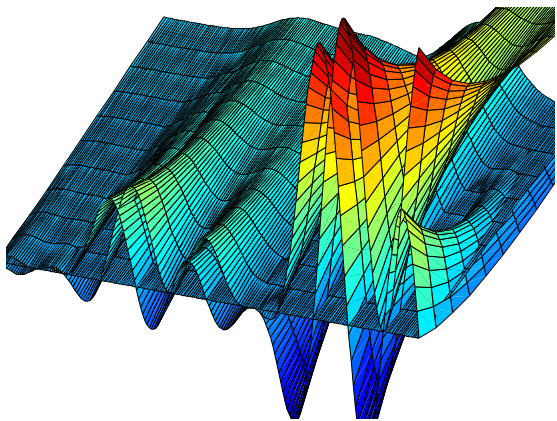}
\caption{$\zeta(\sigma+it)$ for $\sigma \in (1/2,1)$ and $t \in (7001.8,7007.2)$.The surface representing the real component of $\zeta$ is the one that is generally larger.}
\label{6704}
\end{figure}

$\zeta (s)$ is a complex function of the argument $s=\sigma +it$ which has one pole, at $(1+i0)$\cite{derbyshire,edwards,sabbach,vaughn}.  It has ``trivial'' zeros at negative even integers, and at nonuniformly but precisely located positions on the ``critical strip'' $0 \le \sigma \le 1$.  As far as is known, all of these lie on the line $\sigma = 1/2$.  Figure \ref{6704} shows $\zeta$ on the domain $\sigma \in (1/2,1)$ (with $\sigma = 1/2$ on the front edge) and $t \in (7001.8,7007.2)$ (from left to right).  Each point of this figure is determined by an Argand diagram of about 2000 steps if summed directly, or 33 steps if the symmetry investigated here is exploited.  For large $\sigma$, $\zeta$ approaches $1+i0$; in general, the real component is distinguished by being larger than the imaginary component.  ``Flutes'' of increasing amplitude form in both the real and imaginary sheets as $\sigma$ is decreased across the critical strip.  Note that there are several places on the front edge where the real and imaginary components are simultaneously zero.  Two of these near the center are quite closely spaced, too closely to be discerned clearly.  This specific region will be discussed in a later section.  Technical (and popular) interest in $\zeta$ centers on Riemann's Hypothesis (RH) that all of these zeros lie on the line $\sigma = 1/2$.

 The importance of $\zeta(s)$ in the investigation of primes is due to its multiplicative properties.  As can be seen by explicit multiplication, the $nth$ term of the regrouped products of $\zeta^2$ is the number of \textit{divisors} of $n$.  This of course depends on the prime decomposition of $n$, leading to a connection with the primes. As another example, the inverse of $\zeta$ is the Dirichlet series whose coefficients are the \textit{M\"obius} function $\mu (n)$, equal to 1 when n is 1 or a product of an even number of distinct primes, -1 if n is a product of an odd number of distinct primes, and 0 otherwise.

 Riemann exploited an even more dramatic circumstance; the \textit{nth} Dirichlet coefficient of $log(\zeta)$ equals $1$ if $n$ is a prime, $1/k$ if $n$ is the \textit{kth} power of a prime, and $0$ otherwise.  This is easily seen by taking the logarithm of Euler's ``Golden Product'' and expanding the terms;

 \begin{equation}
 \zeta(s) = \Pi_p [1-\frac{1}{p^s}]^{-1}
 \end{equation}

 The coefficients of a nonlinear function of $\zeta$ can be determined using a \textit{Mellin transform}, characteristic of the integral measure $dx/x$, which is invariant under $x \to constant \cdot x$.\footnote{Historically, three nonlinear operations were applied to $\zeta(s)$: $1/\zeta(s)$ (Merten), $log(\zeta(s))$ (Riemann) and $-\zeta '(s)/\zeta (s)$ (Von Mangoldt).  While each of these converts the zeros of $\zeta(s)$ to poles, the properties described above are central to its use in number theory.} Coefficients can be isolated by the orthogonality properties of their characteristic ``log sinusoidal'' variation due to the imaginary component of $s$.  Consequently, transforms of nonlinear functions of $\zeta$, such as $log(\zeta)$ as used by Riemann, show peaks at the prime numbers.  Leaving important details of this standard development to the references\cite[Sec. 1.18]{edwards}, this process leads to a sum of the form

\begin{equation}
f(x)=\sum_{\alpha} \frac{cos(\alpha_i log(x))}{x^{1/2}log(x)} + \ldots
\end{equation}

The coefficients $\alpha_i$ in this transform are the zeros of $\zeta(s)$, the arguments at which $\zeta(1/2+i\alpha_i)=0$.\footnote{The integral of this sum, plus $1/log x$ and some small higher order terms, gives a "staircase" with a \textit{unit rise} of the steps at the primes, and a rise of $1/k$ when $x=p^k$.  For this work, the sum above will be of more immediate interest.\cite[sec.~1.13 to 1.17]{edwards}} The first 20 terms and their sum (ignoring the relatively weak dependence of the denominator $x^{1/2}log(x)$), are shown in Figure \ref{sum20}. The individual terms are ``stacked'', and their sum is shown as the bottom curve.  Note that the minimum at 8 is about 1/3 as deep as the others, and the ones at 4 and 9 are about 1/2 as deep.  This can be regarded as a plot of the coefficients of $log(\zeta(s))$ with an arbitrary vertical scale:

\begin{figure}[pt]
\centering
\includegraphics[angle=0,width=.6\textwidth]{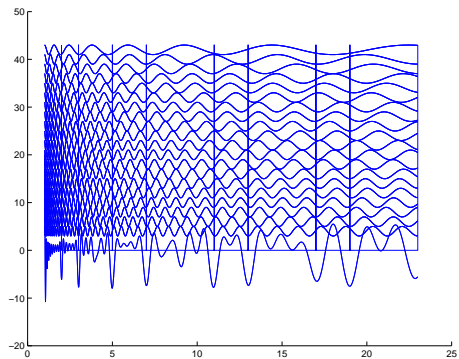}
\caption{$\sum_i cos(\alpha_i log x)$ for the first 20 zeros.}
\label{sum20}
\end{figure}

\begin{figure}[tp]
\centering
\includegraphics[angle=0,width=.6\textwidth]{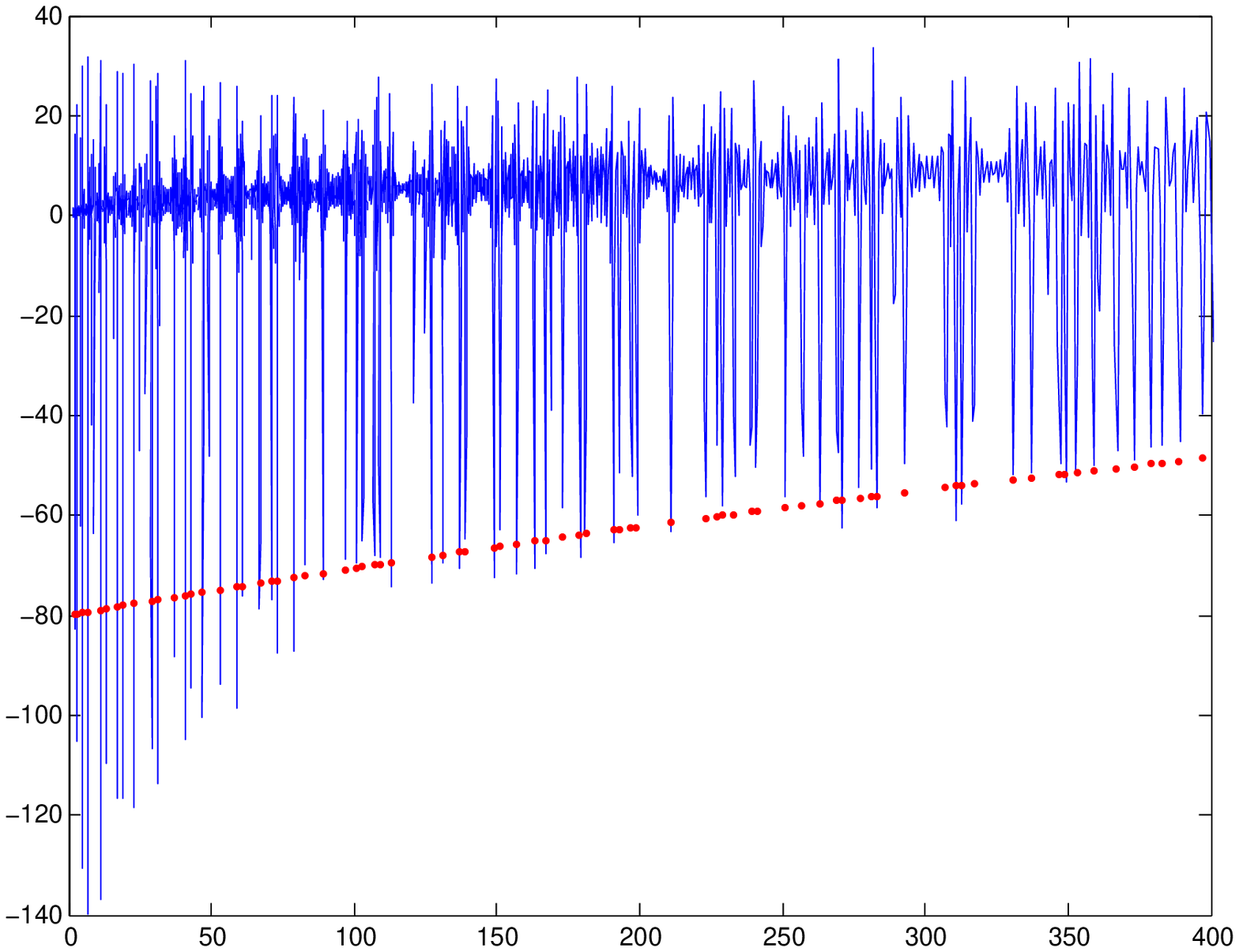}
\caption{$\sum_i cos(\alpha_i logx)$ for the first $10^4$ zeros.}
\label{sum10to4}
\end{figure}

The minima at the primes (and at their $kth$ powers, at reduced height $1/k$) become sharper as more terms are added as shown in Figure \ref{sum10to4}.   Primes are indicated by the red dots.

\begin{figure}[pt]
\centering
\includegraphics[angle=0,width=.6\textwidth]{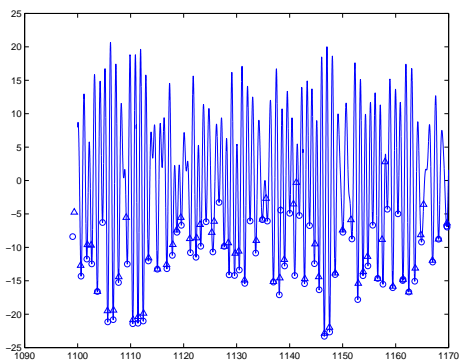}
\caption{$f(y)$ with $\alpha_i \in (1100,1170)$ indicated by ${\mathcal\triangle}$, minima by ${\mathcal\circ}$.}
\label{alphas}
\end{figure}

These figures illustrate quite dramatically the relationship of $\zeta(s)$ to the prime numbers. It is possible to ``invert'' Riemann's expression mathematically, to determine the zeros $\alpha_i$ \textit{given} the prime numbers.  To do this, multiply his equation for $\sum\delta(x-p)$ by $cos(y log x)$ and integrate over x.  This leads to the approximate solution:

 \begin{equation}
 \sum_p^P cos(y log(p)))= \sum_i \psi(\xi_i)
 \end{equation}

 where $\xi_i = y- \alpha_i$, $\psi(\xi) \simeq A sin[\xi_i log P] /\xi_i-sin[\xi_ilog p ]/\xi_i$ and

 \begin{equation}
 A \simeq \frac {1}{log(P/p)} \int_{log(p)}^{log(P)} e^{z/2} \frac {dz}{z}
 \end{equation}

\noindent whose minima approximate the $\zeta$-zeros $\alpha_i$, shown in Figure \ref{alphas}. This method of determining $\alpha_i$ using the primes as ``input data''  provides another visualization of the relationship of the primes to $\zeta$, and could lead to new interpretations of the spacings between the zeros.

\section{Discrete Sampling}

\begin{figure*}[pt]
\centering
\includegraphics[angle=90,width=1\textwidth,height=!]{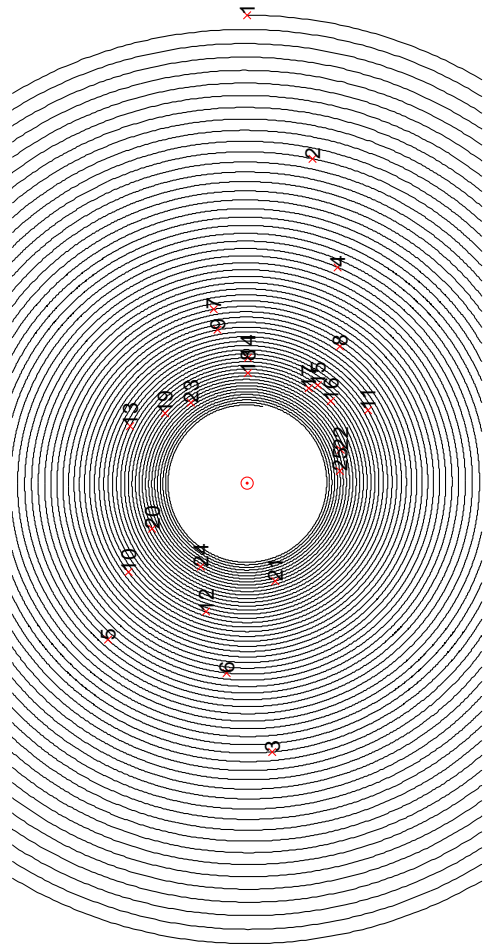}
\caption{Discrete Sampling of $x^{-(1/2+100.586 i)}$ }
\label{discsample}
\end{figure*}

 Line segments, or \textit{steps}, arise from discrete sampling of the smooth function $x^{-s}$.  Step \emph{lengths} depend on the real component of the argument $\sigma$, and vary monotonically as $n^{-\sigma}$.\footnote{For most calculations, $\sigma = 1/2$, and discussions here are generally restricted to the ``critical strip'', $\sigma \in (0,1)$.}  Step \emph{angles}, however, depend on the imaginary component of the argument \textit{t} through the relation $\theta_n = -t log(n)$.  This involves modular arithmetic to reduce the angle to an interval of $2\pi$.  Because  \textit{t} can be quite large, $\theta_n$ generally exceeds any possible domain of convergence of a conventional Taylor series.\footnote{In addition to the usual limitations of computer speed and human patience, numerical evaluation of the modulus function requires enough digital capacity to support $t log(n)$ mod $2\pi$ accurately for large \textit{t log(n)}, since it is essentially the difference of two large numbers.  Numerical techniques are available to increase the range, if desired.\cite[sec.~7.1]{recipes}}  Discrete sampling of the continuous spiral $x^{-s}$ at $n \in N$ is illustrated by Figure \ref{discsample} for $s = .5 + 100.586 i$, a value simply chosen to demonstrate key points without overlapping step labels.  In spite of the large number of "skipped" cycles in the initial steps of the spiral, there are accurate difference equations relating successive angles.  The regularities in these discrete equations responsible for the observed symmetries in the Argand diagram arise in much the same manner as discrete sampling of the time dependent equation $dP/dt = k P (1-P)$ leads to the intricacies of the chaotic logistics map \cite[sec.~1.1]{devaney}.

The \textit{nth} step is the vector from the origin to the numbered point; for the $\zeta$ function the first step is always 1. In this example, the angle $\theta_2 = -t log(2)$ is -69.7211 radians, slightly more than 11 complete circles. The mod function reduces this value to -.6061 radians.  As n increases, the number of complete circles between adjacent steps decreases, becoming less than 1 when $n > t/{2\pi}$, about 16 in this example.  Note that steps 15, 16, and 17 are almost in the same direction, having angle changes of approximately $2\pi$.  Also, as n increases beyond 17, the differences in angle become quite regular. It may seem that the reduced angle differences for $n \ll t$ resemble random number generator schemes that keep only the least significant digits.  This process is far from random, however; it will be shown that $\theta_2$ is equal to the total angle change for an appropriately chosen range of ``conjugate steps'' of larger \textit{n}.

\begin{figure}[pt]
\centering
\includegraphics[angle=0,width=.6\textwidth]{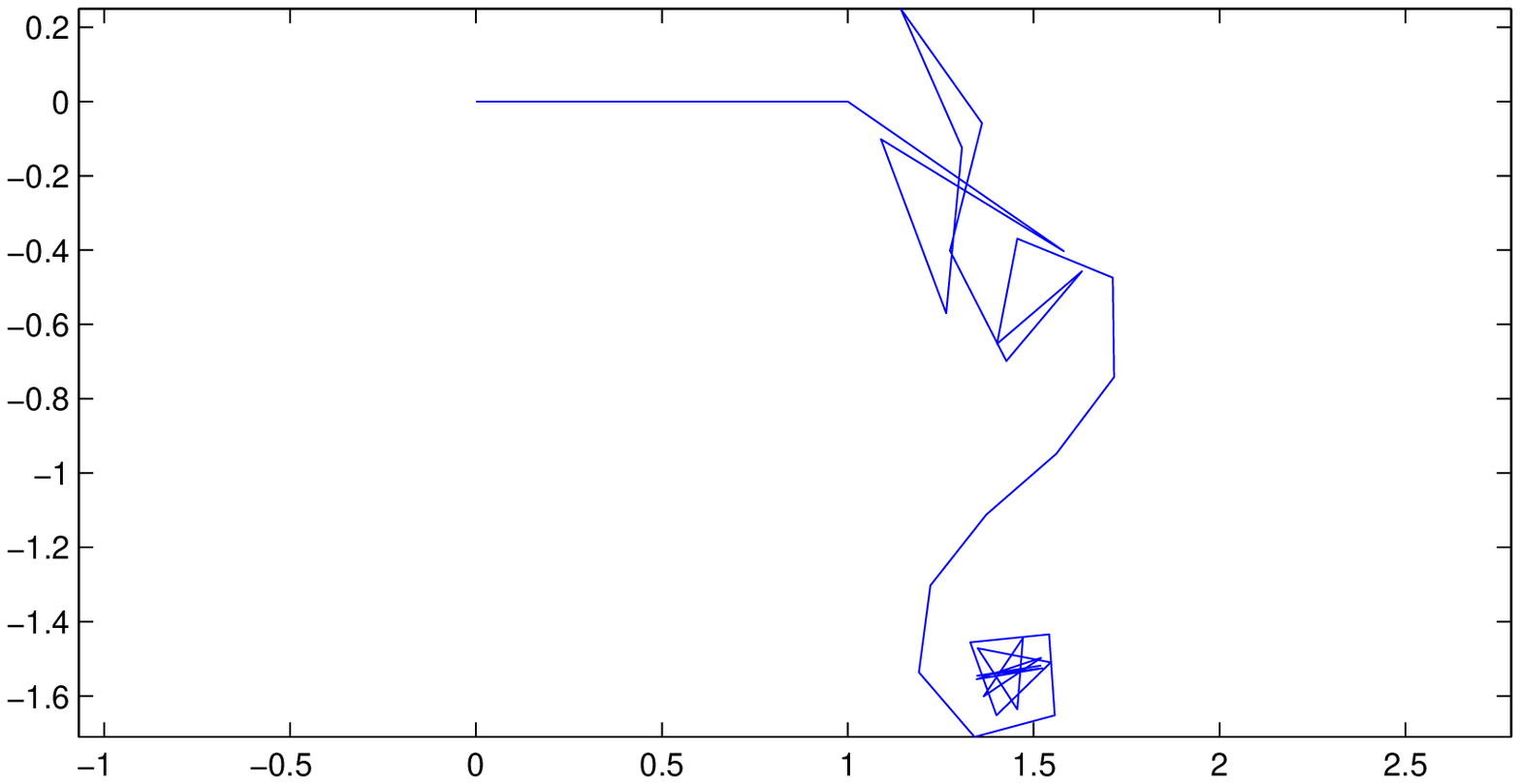}
\caption{Argand Diagram for s=1/2+100.586 i }
\label{argandsmpl}
\end{figure}

Figure \ref{argandsmpl} is the Argand diagram for this discrete sampling example.  While suppressing the ``$2\pi$ jumps'' of the modular arithmetic, it clearly displays the reversals such as at steps 10 and 11, and the smooth arc of steps 15, 16 and 17.  The structure from step 10 to the ``star'' near $(1.5,-1.5)$ is a prominent recurring feature, and the ``converged'' value of $\zeta (1/2+100.586)$ will be shown to be the geometric center of the final star.

\section{Step Angles}

The general features for $t \simeq 100$ are maintained as $t$ increases to very large values, but many characteristic structures arise.  It is useful to have mathematical descriptions of the angle differences between steps. As mentioned above, a Taylor series expansion of $x^{-s}$ would be valid over a much smaller domain than is required.   Instead, the discrete angles and their differences are given by the direct use of logarithms and their small argument expansions.  This yields the ``Discrete Taylor Series'' (DTS) whose coefficients are discrete angle differences rather than derivatives:

\begin{equation}
\theta = -t log(n)
\end{equation}
\begin{equation}
\delta \theta = -t log(\frac {n+1}{n}) \simeq -\frac{t}{(n+1/2)} \end{equation}
\begin{equation}
\delta^2 \theta  \simeq \frac {t}{(n+1/2)^2}
\end{equation}
\begin{equation}
\delta^3 \theta \simeq - \frac{2 t}{n^3}
\end{equation}

Note that these are \emph{forward differences}, in accordance with the digital sampling definition for $\zeta (s)$: the angle of the $nth$ step is determined by the direction from the origin to the spiral at the \emph{beginning} of the step. The general pattern to advance each term is

 \begin{equation}
 \delta\theta^{k}_{n+1}=\delta\theta^k_{n} + \delta\theta^{k+1}_{n} \qquad k \in Z
 \end{equation}
\noindent
 where $k=0$ corresponds the step angle and $k > 0$ represents the \textit{kth} angle difference.  Each coefficient can be considered as having been reduced modulo $2\pi$.  With these definitions of the differences, the result of a succession of such steps centered at step $n_0$ is the  DTS.\footnote{These "binomial-like" coefficients are due to the definition of the forward differences.  They are derived by an inductive application of the difference algorithm pattern. For centered differences, the usual Taylor series coefficients result.}

\begin{multline}
\theta (n_0 + \delta n) = \theta (n_0)+ \delta n \delta \theta (n_0)
+\frac {1}{2}\delta n(\delta n -1)\delta ^2 \theta (n_0) \\
\qquad +\frac {1}{6}\delta n(\delta n -1)(\delta n -2)\delta ^3 \theta (n_0)+\ldots
\end{multline}

   Because a large number of $2\pi$ cycles may be ignored in the modulus operation, it may seem that this representation could not possibly give accurate angle differences for the discrete sampling procedure.  To demonstrate numerically that it does, values of the DTS expansion about the point $n_p$ for $t=10^9$ are shown in Table \ref{dtstable}.

\begin{table*}[hp]
\begin{tabular}{|c|c|c|c|c|}
\hline
Step & $\theta$ & $\delta \theta$ & $\delta^2 \theta$ & $\theta (series)$ \\
\hline
12611 & $-(n_{12615}+50474)\times2 \pi - 4.6057$ & $-(n_p +4)\times2\pi - 5.1937$ & $2 \pi +.0036$ & $13\times2\pi -4.6057 $\\
12612 & $-(n_{12615}+37854)\times2 \pi - 3.5162$ & $-(n_p +3)\times2\pi - 5.1900$ & $ 2\pi +.0027$ & $8\times2\pi -3.5162 $\\
12613 & $-(n_{12615}+25235)\times2 \pi - 2.4230$ & $-(n_p +2)\times2\pi - 5.1874$ & $2 \pi +.0017$ & $4 \times2\pi -2.4230$ \\
12614 & $-(n_{12615}+12617)\times2 \pi - 1.3272$ & $-(n_p +1)\times2\pi - 5.1857$ & $2 \pi +.0007$ & $ -1.3272$ \\
12615 & $-1502843128\times2 \pi - 0.2297$ & $-12615\times2\pi  -5.1851$ & $2 \pi -.0004$ & $-0.2297$ \\
12616 & $-(n_{12615}-12615)\times2 \pi - 5.4148$ & $-(n_p -1)\times2\pi - 5.1854$ & $2 \pi -.0023$ & $-5.4148$ \\
12617 & $-(n_{12615}-25230)\times2 \pi - 4.3170$ & $-(n_p -2)\times2\pi - 5.1867$ & $2 \pi -.0033$ & $ -4.3170$ \\
12618 & $-(n_{12615}-37844)\times2 \pi - 3.2206$ & $-(n_p -3)\times2\pi - 5.1891$ & $2 \pi +.0036$ & $ -3.2206$ \\
12619 & $-(n_{12615}-50457)\times2 \pi - 2.1265$ & $-(n_p -4)\times2\pi - 5.1924$ & $2 \pi +.0043$ & $2\times2\pi -2.1265$ \\
\hline
\end{tabular}
\caption{Angles and angle differences near step  $n_p = 12615$, including multiples of $2\pi$ that are removed by $N mod,2\pi$.}
\label{dtstable}
\end{table*}

This table also makes it apparent why 4 places of accuracy in the step angle requires an arithmetic register with a capacity of at least 15 places for $t=10^9$.  For \emph{each} of these steps, $\delta \theta$ encompasses a total angle change of $n_p$ complete circles.\footnote{Since the DTS will generally be used only for $n>n_p$, this is a "worst case" example. The accuracy is much higher when the step number of the expansion center is larger, because the $kth$ angle difference coefficient decreases as $n^{-k}$.}

\begin{figure}[pt]
\centering
\includegraphics[angle=0,width=.6\textwidth]{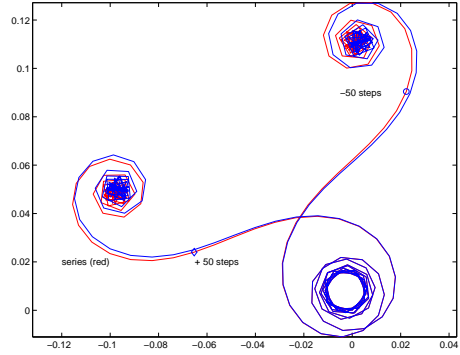}
\caption{Comparison of actual steps and series near $n_p$}\label{pendseries}
\end{figure}

A graphical example of this is presented in Fig. \ref{pendseries}, an extension of the range detailed in Table \ref{dtstable}.  The angles determined by the DTS taken to $\delta ^4\theta$ for this case, are shown in red.  Note the eventual asymmetry in the arc distance to steps $\delta n = \pm 50$, due to the shortening of steps as $n$ increases.

The $\delta \theta$ relationships can be inverted to determine the step number at which particular angle differences occur.  Expressions for these are:\footnote{The angles used in these equations are \emph{not} to be reduced modulo $2\pi)$; for large \textit{n} higher angle differences may be $< 2 \pi$, obviating the need for the modulus function.  When integer step values are required, use the \textit{fix} function $[n]$, and correct with the remainder $\delta n = n-[n]$ as necessary.}

\begin{equation}
n(\theta) = e^{-\theta/t}
\end{equation}
\begin{equation}
n(\delta \theta) = \frac{1}{e^{-\delta\theta/t}-1}
\end{equation}
\begin{equation}
n(\delta ^2\theta)\rightarrow n^2+2n=\frac{1}{e^{\delta^2\theta/t} -1}
\end{equation}
\begin{equation}
n_{pk}(\delta ^2\theta = 2 \pi k) \simeq \sqrt {\frac {t}{2 \pi k} } + \ldots
\end{equation}

Referring to Figure \ref{discsample}, the step number for $\delta \theta = 2 \pi$ is 16.  At this step number, the second difference in angle is small, leading to the close spacing of steps 15, 16 and 17.  Steps greater than 16 have a first difference less than $2\pi$, and steps greater than $n_p=\sqrt{t/(2\pi)}$ have second difference less than $2\pi$.  This is the center of symmetry, a particularly important step number for this discussion.

\section{Enumeration of the Symmetry Features}

Almost immediately, initial steps $n$ and opposing symmetric steps $\tilde{n}$ are seen to be related by a reciprocal relation of the form $n \tilde{n} \simeq t/2 \pi$.  To the next order of approximation,

\begin{equation}
(n+1/2)(\tilde{n}+1/2) = \frac{t}{2\pi}
\end{equation}
\noindent
so that the origin is reciprocal to $[t/\pi] - 1/2$ and step $[(t/2\pi)^{1/2}]-1/2$ is reciprocal to itself.

The key features of the symmetry will now be systematically discussed and illustrated.  These are ordered by the increasing "multiplicity" M of conjugate to initial steps, as defined below.  Each is demonstrated graphically within the Argand diagram.

\begin{enumerate}
\it
\item The center of the symmetry is a ``pendant'' near step number

    \begin{equation}
    n_p = [\sqrt{t/(2\pi)}],\quad at \quad which \quad \delta^2 \theta = 2 \pi
    \end{equation}
\rm
In Figure \ref{xy9}, this feature is within the central blob near (-0.9,-1.650).  It is more easily seen in the center of Figure \ref{blob1pt6}. Although this quantity appears in the analytical literature and may well have an established name, it is referred to in this context as $n_p$ for ``\emph{pendant}'' because of its appearance when $t$ is large.  The pendant holds a ``star sapphire'' if $\delta \theta (n_p)$ is close to an odd multiple of $\pi$, or a ``pearl'' if it is close to an even multiple of $\pi$.  For $t=10^9$, $\delta \theta (n_p) \simeq \pi/3$, giving a nested set of polygons resembling a faceted diamond (see Figure \ref{pendseries}). Algorithms to determine $\delta \theta (n_p)$ and the precise center of symmetry in terms of $p=\sqrt{t/2\pi} - [\sqrt{t/2\pi}]$ are presented below\footnote{Special effort is required, however, when $\delta \theta (n_p) = 0$.}.
\\

\it
\item Each single step $n < n_p$ is ``mirrored'' by $M \simeq  n_p^2/ n^2$ conjugate steps $\tilde{n} > n_p$ where:

\begin{equation}
 \frac{n_p^2}{n +1/2} \le \tilde{n} \le  \frac{n_p^2}{n -1/2}
 \end{equation}

\rm
 This is a fundamental element of the symmetry, and it will be discussed in detail from many standpoints. Figure \ref{xy9} shows the effect clearly, most convincingly by comparing initial steps and final spirals.  M varies greatly across the diagram, from 1 at $n_p$ to $2t/3\pi$.  The region conjugate to the first step contains 2/3 of \emph{all} the steps shown in the diagram.

\it
\item When $M \simeq 1$, the initial and conjugate steps are of nearly the same length for all
 $\sigma$, with step angles given by an odd function of $n-n_p$.

\rm

 Near $n_p$, the second difference DTS term with coefficient $\delta n (\delta n -1)/2$ always gives a multiple of $2 \pi$.  Thus, the the first nonzero \emph{even} function of $\delta n$ is the term $-1/2 (\delta_n)^2$ in the coefficient of $\delta^3 \theta$.  For small $\delta n$ near $n_p$, the steps lengths do not vary appreciably, leading to the central symmetry.

Consider the original step $(n_p-1)^{-s}$ and its neighboring conjugate step $(n_p)^{-s}$.  One of the most fundamental aspects of the $\zeta$ function relates to the \textit{inversion} $s \to 1-s$. When $\sigma = 1/2$, this simply gives the complex conjugate of the argument $s$. Under this transformation, the leading term of the ratio of the conjugate step to the transformed original step is $n_p^{1-2s}$.  It will be shown below that this relation between conjugate steps and transformed original steps, with a carefully defined complex angle, is independent of M.

\it
\item For $M > 1$, the ``vector sum'' of the conjugate steps equals the length of the initial step if $\sigma = 1/2$. Each initial step and its conjugate step(s) are bilaterally symmetric with respect to angle about an axis through the center at $n_p$, lying in a direction perpendicular to $\Theta = -t log n_p+t/2+\pi/8+\ldots,( \mod 2\pi) $.

\rm
   The first term in the expression for the symmetry angle is just twice the angle of the step at $n_p$.  The additional terms arise from discrete effects when $n_p \notin P$ and the "tilt" of conjugate regions, as will be discussed below.

\it
\item The Vector Sum
\rm

   The angle relation is seen by multiplication of the logarithm of the reciprocal relation by $-it$.  There is a scaling principle that demonstrates the ``vector sum equality''.    Consider an original step at $n_p/M^{1/2}$, of length $M^{\sigma/2}n_p^{-\sigma}$. Its conjugate region at $M^{1/2} n_p$ contains M steps, each of length $(M)^{-\sigma/2} n_p^{-\sigma}$.  Vector addition of these M steps in the complex plane generally averages to $M^{1/2}$ steps.

  This can be demonstrated explicitly using trigonometry when \textit{M} is a small integer, or numerically for larger $M$.  For $M \gg 1$, it can be approximated using the DTS. It is also similar to the result of the traditional argument for a random walk, in which $[\sum_1^M \overrightarrow{\lambda_i}]^2 \simeq M \lambda^2$ with vector steps $\overrightarrow{\lambda}$ having uncorrelated directions or the uniformly spaced roots of unity. Of course, in the present application, the direction and distance of the summed terms have well-defined values.  A precise expression will be derived below, using ``classical'' analytic number theory.  The result of this discussion is that the vector sum of the conjugate steps has a length $M^{(1/2-\sigma)/2}n_p^{-\sigma}$, which is equal to $n_p^{-\sigma}$ only for $\sigma = 1/2$.

   For an original step at $n_p/M^{1/2}$, the length of the conjugate region is proportional to the initial step length under the substitution $\sigma \to 1-\sigma$, which is

    \begin{equation}
    M^{1/2}M^{-\sigma/2}n_p^{\sigma -1}
     \end{equation}

     This was the result for adjacent steps when $M = 1$, and is now seen to apply for larger M.  When the angles due to $it$ are included as well the lengths due to $\sigma$;

  \begin{equation}
  \sum_{\tilde{n}}( n')^{-s} = n_p^{1-2s} n^{s-1} \qquad
   \tilde{n} \in (\frac{n_p^2}{n+1/2},\frac{n_p^2}{n-1/2})
  \end{equation}

  When $\sigma=1/2$, this has the result that the amplitudes of the conjugate region and original step are equal, and the sum of their phases is $e^{-it log(n_p)}$.\footnote{The traditional mathematical term for the magnitude of a complex number is ``modulus''.  Since that term is used here in another context, the expression generally used in waveform analysis, $amplitude \cdot e^{i \cdot phase}$, will be used for the polar representation of complex numbers.}  This relation holds for \textit{each} original step $n \le n_p$ and its conjugate region.  It is emphasized that this is a \textit{detailed} expression, valid for each of the $n_p$ steps leading up to the symmetry center, not only a global result for a sum over the entire range of steps.

\it
\item Rational Fraction Conjugacy
\rm

It is not always the case that a \emph{single} original step corresponds to an integer number of conjugate steps. Interesting features arise when M is a rational fraction and $k$ original steps are conjugate to $Mk$ steps.  For example, Figure \ref{saw45} compares two opposing "circular saw blades" having 4 and 5 steps per "tooth".  Figure \ref{blob1pt6} shows the entire region around $n_p$ with $M \le 1.6$ for $s=1/2+10^9$. The center of symmetry and many "rational fraction" conjugate pairs are evident in this rather impressive figure.  Note that  conjugate features having a different number of steps than the original region have the same scale size, as expected for $\sigma = 1/2$.

\begin{figure}[pt]
\centering
\includegraphics[angle=0,width=.6\textwidth]{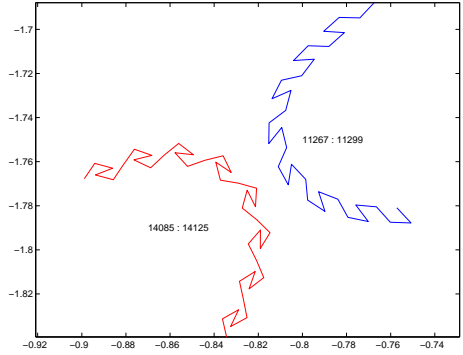}
\caption{Opposing features for $M=5/4$}
\label{saw45}
\end{figure}

\begin{figure}[pt]
\centering
\includegraphics[angle=0,width=.9\textwidth]{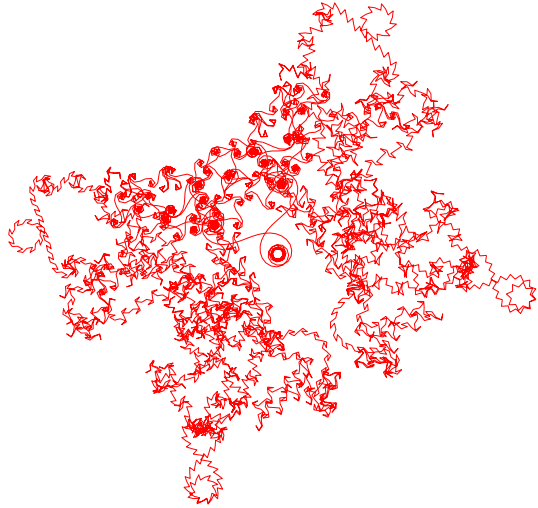}
\caption{The region with $M < 1.6$ for $s = 1/2+10^9 i$}
\label{blob1pt6}
\end{figure}

\it
\item Small Integer Conjugacy
\rm

 Consider $M$ a small integer.  The original step at $n_p / M^{1/2}$ has angle $\theta_n=\theta_p -t/2 log M$, where $\theta_p \simeq -t log n_p)$. \footnote{Because of the extreme sensitivity of this important angle, it will be defined more precisely below.}.  Using $M = 5$ as a simple example where the conjugate step angles can be added explicitly, each original step transforms to the 5 intervals with the successive angles

\begin{multline}
  \theta_0 \Rightarrow \theta_0+1/5\delta \theta \\ \Rightarrow \theta_0+2/5\delta \theta +1/5(2\pi) \\
  \Rightarrow \theta_0+3/5\delta \theta +3/5(2\pi) \\
  \Rightarrow \theta_0+4/5\delta \theta +1/5(2\pi)
\end{multline}

Here, $\theta_0$ is the angle of the first conjugate step, and $\delta\theta$ is the step angle change at the \textit{original} step. The first difference terms add to the original step angle change, and the second difference terms add to $2\pi$.  These five steps constitute the region conjugate to the original step, as shown in Figure \ref{fivesteps}. (The characteristic ``Z'' pattern of these 5 steps is also evident in any region where M is a rational fraction involving 5.)

\begin{figure}[pt]
\centering
\includegraphics[angle=0,width=.6\textwidth]{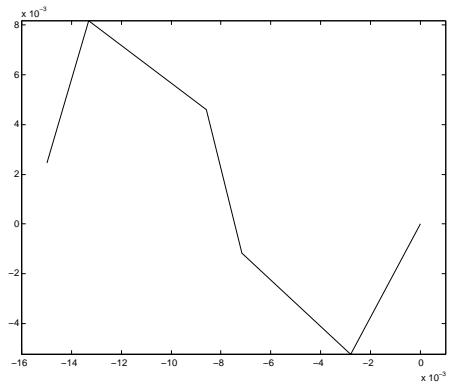}
\caption{The 5 Steps Conjugate to $n=n_p/\sqrt{5}$ for $s=1/2+10^9 i$}
\label{fivesteps}
\end{figure}

Because step $n_p/\sqrt{5}$ is also a point where $\delta\theta^2 = 10 \pi \to 0  mod 2\pi$, the steps surrounding this point give a copy of the pendant at $n_p$.  This feature and its conjugate region are seen in Figure \ref{xy9} at $(-.85,-1.01)$ and $(-1.48,-1.42)$ respectively.  This "minipendant" and its its conjugate have been enlarged and translated into proximity for comparison in Figure \ref{pair5}.  The general results of the "conjugate step multiplication" are easily seen.  The closed loop consists of 48 original steps and 240 conjugate steps.  It will be shown below that the distance between points representing the the centers of the vestigial spirals at each end of the ``Z'' equals the original step length.

\begin{figure}[pt]
\centering
\includegraphics[angle=0,width=.6\textwidth]{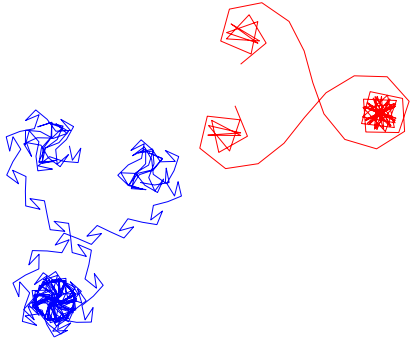}
\caption{Conjugate "minipendants" for $M=5$}
\label{pair5}
\end{figure}

For $s=1/2+10^9 i$, there happen to be two conjugate features that almost exactly meet on the symmetry axis without translation, as shown in Figure \ref{kiss14}.  The ratio of the central step numbers shown in the figure indicates that $M=13.8$.  This is consistent with a visual estimate of the conjugate steps, and is one of the clearest examples of the conjugacy and symmetry of the Argand diagram.

\begin{figure}[pt]
\centering
\includegraphics[angle=0,width=.6\textwidth]{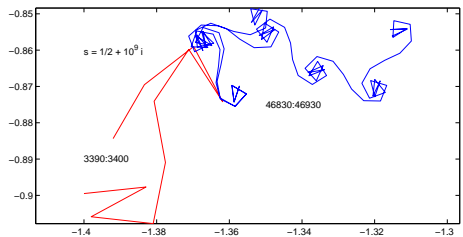}
\caption{Opposing features for $M \simeq 13.8$}
\label{kiss14}
\end{figure}

The examples of rational fraction and small integer values for M explain the apparent lack of symmetry for original steps and their conjugate regions.  For $\sigma = 1/2$ the conjugate regions have a total extent and net angle corresponding to the original step.  This is demonstrated for large M in the next item.

\it
\item For $ M \gg 1$, the conjugate steps form a pair of opposing spirals linked at the center by
an inflection point where adjacent step angles differ by an even multiple of $\pi$. This generic shape begins to be evident for $M \ge 5$, and is the most distinctive feature for larger values of M.

\rm

 This is a dominant feature of Figure \ref{xy9}, and is clearly explained by the DTS. Defining the center of the \textit{kth} conjugate region as the step $\tilde{n}_k$ at which the first angle difference is $2\pi k$ gives $\tilde{n}_k \simeq t/2 \pi (k- 1/2)$. Each conjugate region extends to the nearest adjacent steps at which $\delta \theta$ is an odd multiple of $\pi$, or

 \begin{equation}
  \tilde{n}-\frac{t}{(2k +1)\pi} \to \tilde{n}+\frac {t}{(2k -1)\pi}
 \end{equation}

    Using the second angle difference at $\tilde{n}$ of $2 \pi k^2/n_p ^2$, the DTS shows that the first angle differences approach $\pm \pi$ at the extremes of the conjugate range. It will be shown that the angle between the central step of a conjugate region and the ``center-to-center'' extent of the region is $\pi /4(1-k^2/n_p^2)$. Note that as \textit{k} increases from $1 \rightarrow n_p$ there are $n_p$ regions conjugate to each single original step $k$.\footnote{We refer to this conjugate structure loosely as an Euler spiral, in which the step angle varies as the square of the arc distance.
    These steps approximate that geometric form when the third angle difference in the DTS is small and the variation of step length with n is ignored\cite[pp. 1-3]{orangepeel}. }

While distances require explicit vector summation of steps, angle differences are easily defined over many steps.  For example, use of Equation 11 shows that the angle between the steps at the inflection points of adjacent Euler spirals $k$ and $k+1$ is

\begin{equation}
t log[ e^{-2\pi k}-1]-t log [ e^{-2\pi (k+1)}-1] \simeq -t log(1+\frac{1}{k})
\end{equation}
\noindent
which equals the first angle difference of the corresponding original steps.  This is yet another visualization of the detailed symmetry.

\it
\item The center of the \textit{kth} scroll, where $\delta \theta = (2k-1)\pi$, is given to first order by

 \begin{equation}
 \sum_1^{N_k} n^{-s} - N_k^{-s}/2 +\frac{\sigma + i \delta t}{4 N_k^{1+\sigma}}
 \end{equation}
\noindent
 where $N_k = [t/[(2k-1)\pi+1/2]$ and $\delta t = t-(2k-1)\pi N_k$.

\rm
The precision of the symmetry may not be visually evident, because of the central tilt angle of $\pi /4$ and the diameter of the spiral scrolls compared to the original single step.  In order to demonstrate it numerically, it is necessary to develop a geometrical algorithm to locate the center of the $kth$ spiral.   The ``zero-order'' spiral center position is computed by summing the steps to $N = [t/(2k-1)\pi]$ and subtracting $N^{-s}/2$, a procedure often encountered in the traditional analytical development.

 However, steps cross the spiral center in a sparse pattern and the exact center will not usually lie directly on the step. To locate the \textit{kth} spiral center to first order:

 \begin{enumerate}
 \item Determine $ N=[t/(2k-1)\pi]$ and $ \delta t = t-(2k-1)\pi N $.
 \item
     Calculate $ \sum_1^N n^{-s} -N^{-s}/2 +(\lambda_N(\sigma ) + i \tau_N(\delta t)) $.
 \end{enumerate}

 The \textit{first-order} correction has two components; \textit{longitudinal} $\lambda_N$, in the the step direction and \textit{transverse} $\tau_N$, perpendicular to it.  Both depend on $s$, but the arguments indicated stress the dominant dependence.  The longitudinal correction is due to the decrease of the length of steps with $N$ when $\sigma > 0$.  Denoting the ``true center'' of a step by its midpoint $N^{-s}/2 + \lambda_N$ and assuming that $\lambda_N$ varies slowly with $N$, the condition is

 \begin{equation}
 \sum^N n^{-s} - N^{-s}/2 +\lambda_N = \sum^N n^{-s} + (N+1)^{-s}/2 -\lambda_N
 \end{equation}
\noindent
 where $\theta(N) \simeq \theta(N+1) \pm \pi$.  The solution is

 \begin{equation}
 \lambda_N = \frac{\sigma}{4 N^{1+\sigma}}
 \end{equation}

  The transverse correction is assumed to vary linearly with $\delta t$.  If $\delta t$ were 0, the \textit{centered} first angle difference $-t/N$ would be $(2k-1)\pi \simeq 0$.  Using the centered form of the DTS,

  \begin{equation}
  \theta(\delta n) = \theta_0 +\delta n \delta \theta _0 +1/2 \delta n^2 \delta^2_0 \theta + \ldots
  \end{equation}
\noindent
  where here $\delta \theta_0 = -t/N \simeq \pi$ and $\delta^2 \theta_0 = t/N^2$.
   Thus, the preceding and succeeding steps each differ from step N by the identical centered second angle difference $(2k-1)\pi/2N$.  The center position lies \textit{on} this symmetrically located central step.\footnote{Intuitively,it may seem that the spiral center should lie on the two steps that reverse to overlay exactly, with a \textit{forward} first angle difference of an odd multiple of $\pi$. Some reflection will show that it is the step whose \textit{centered} first angle difference is an odd multiple of $\pi$ that contains the center.}  When $\delta t=(2k-1)\pi$ (which is equivalent to $\delta n = 1$) the center lies on the midpoint of the next step, $N+1$.  It is therefore displaced laterally by $[(N^{-s}/2)((2k-1)\pi/2N)$.  This gives;

   \begin{equation}
   \tau_N = \frac{i \delta t}{4N^{1+s}}
   \end{equation}

 To test this result numerically, spiral centers were calculated to first order for Figure \ref{xy9} with $s=1/2+10^9i$.  Defining $Q_k(s)$ as the distance between adjacent spiral centers for the $kth$ Euler spiral, divided by the original step under the inversion $s \to 1-s$, gives the amplitudes and phases shown in Figure \ref{ratiopt5}, down to approximately $M=5$.  One sees that $Q$ is a universal ratio, independent of $k$ (or $M$)\footnote{Note that, since the first step is of unit length, the extent of the final Euler spiral equals $Q(s)$.}.

\begin{figure}[pt]
\centering
\includegraphics[angle=0,width=.6\textwidth]{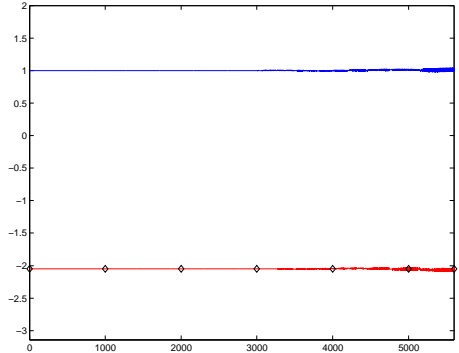}
\caption{Complex ratios of conjugate Euler spirals to original steps, $s=1/2+10^9i$}
\label{ratiopt5}
\end{figure}

For $\sigma = 1/2$, the length ratios are unity to within $10^{-6}$ for small k, rising to an error of a few percent at the end of the range with the first order correction derived above. The sum of conjugate and original angles is a constant whose analytic value, which will be discussed below, is indicated by the red $\diamond$s.

For $\sigma = .45$, a similar calculation illustrates the ratio of the conjugate region to the original step under the inversion $s \to 1-s$.  Because the angles depend only on $t$, they the same as for $\sigma = .5$.  The amplitude ratios in Figure \ref{ratiopt45} reproduce the factor $n_p^{0.1}$, shown by $\diamond$s, to seven or more places for small k.  When $M=5$, even that vestigial Euler spiral is accurate to a few percent with just this first order term.

\begin{figure}[pt]
\centering
\includegraphics[angle=0,width=.6\textwidth]{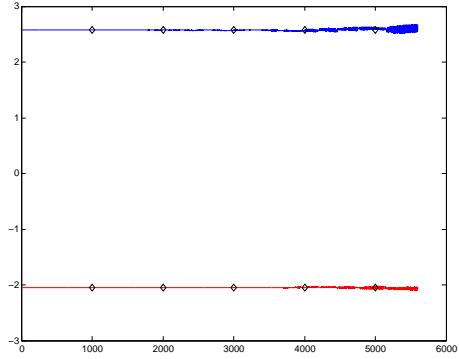}
\caption{Complex ratios of conjugate Euler spirals to original steps, $s=.45+10^9i$}
\label{ratiopt45}
\end{figure} \

\it
\item The center of the pendant, denoted here by $P(s)$, is given to lowest order by $\sum_1^{n_p} n^{-s}$.  This leads to an error of $\mathcal{O}(n_p^{-1/4})$. The first order correction is

\begin{equation}
P(s) = \sum_1^{n_p} n^{-s} -\frac{e^{-i(t log(n_p)+2 \pi p)}}{2 n_p^{1/2} cos(2 \pi p)}
\end{equation}

\noindent where $p = \sqrt{t/2\pi} -[\sqrt{t/2\pi}]$.
 \rm

  This point is the actual center of symmetry of the entire Argand diagram, and is crucial to understanding the Riemann-Siegel equation below.  When $t$ is extremely large, as in Figure \ref{pendseries}, the pendant center $P(s)$ is easily recognized.  For smaller values of $t$ however, it is less apparent and locating it requires more effort than did the spiral centers, where $\delta \theta$ is always near $\pm \pi$.  In the pendant, various values of $\delta \theta$ are encountered, depending on $p$.  By expressing $\delta \theta$ in terms of $n_p$ and $p$ instead of $t$, multiples of $2\pi$ are easily omitted;

 \begin{multline}
\delta \theta (n_p)= -t log \frac{n_p +1}{n_p} = -2\pi (n_p +p)^2 (1+\frac{1}{n_p}-\frac{1}{2n_p^2}+\ldots) \\
\simeq 2 \pi (1/2-n_p -2 p +\frac{p^2}{n_p}) \to \pi - 4 \pi p
 \end{multline}

 The term $2 \pi n_p$ is zero mod $2\pi$, and $p^2/n_p$ is ignored for reasonably large $n_p$.

  $\delta \theta$ advances through a range of $4\pi$ as $p$ increases from $0$ to $1$.  For $p = 0,1/2$, $\delta \theta = \pi$, representing reversals of steps, while $p=1/4,3/4$ give $\delta \theta = 0$.

 For each step $n$, a point $P_n$ can be defined as the intersection of the lines which bisect the angles at each end of the step. This point is the center of a local osculating circle to step $n$ and its neighbors. In regions of the Argand diagram where $\delta^2 \theta, mod 2\pi$ is not small, these points are widely distributed. For $\delta^2 \theta,  mod 2\pi \simeq 0$, and assuming that step lengths change slowly, adjacent center points congregate closely. For $\delta^2 \theta = 2 \pi$, $\delta \theta$ is the same at each end of the step and $P_n$ is the vertex of an isosceles triangle whose base is the step $n_p$ and whose vertex angle is $\delta \theta$.  The length $L$ is given by $L sin(\delta \theta /2) = n^{-\sigma}/2$, and the direction of the side in contact with the end of step $n_p$ is normal to that of the step, minus $\delta \theta /2$.  The common vertex of these neighboring triangles \textit{is} the pendant center.  As above, it is convenient to express the various angles in terms of $p$ and $n_p$, which gives

 \begin{equation}
P_{n_p} \equiv P(s) = \sum_1^{n_p} n^{-s}-\frac{e^{-i(t log(n_p)+2 \pi p)}}{2 n_p^{\sigma} cos(2 \pi p)}
 \end{equation}

\it
\item\ Higher order terms.
\rm

Algorithms to determine the pendant center require some care.   The array of ``pendants'' in Figure \ref{pend3x3} demonstrates this point for a set of values of $\delta \theta (n_p)$ where $n_p=100$.

\begin{figure}[pt]
\centering
\includegraphics[angle=0,width=1.3\textwidth]{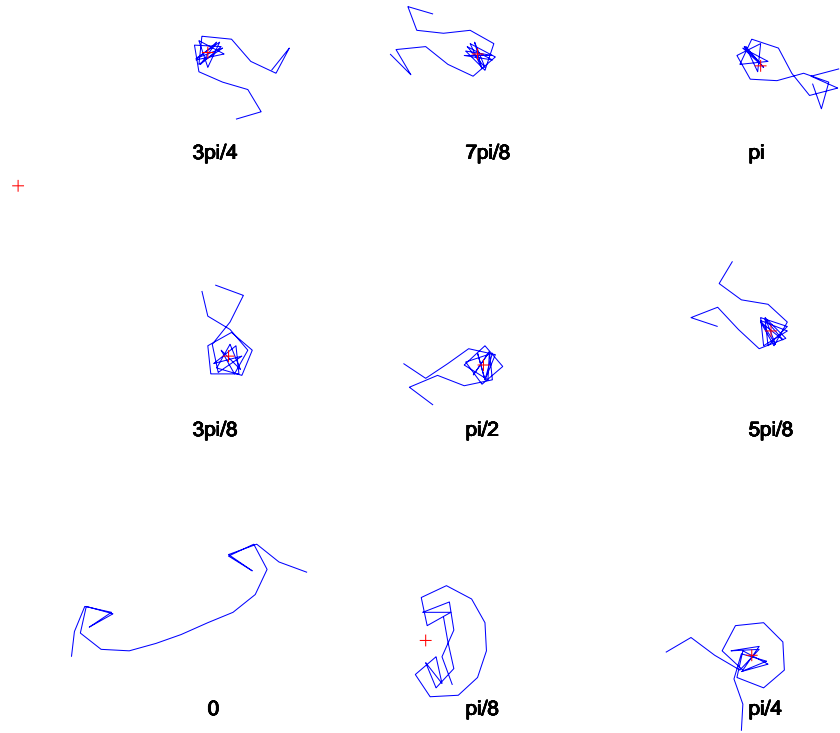}
\caption{Pendant centers near $n_p$ for various values of $\delta \theta$ in steps of $\pi /8$.}
\label{pend3x3}
\end{figure}

 One construction is illustrated in Figure \ref{construct}.

\begin{figure}[tp]
\centering
\includegraphics[angle=0,width=.6\textwidth]{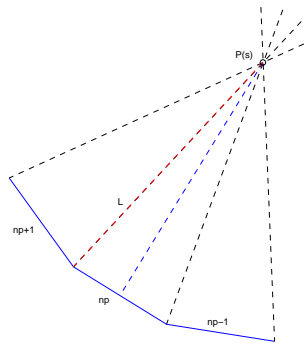}
\caption{Geometrical construction of a pendant center.}
\label{construct}
\end{figure}

If the actual angle differences at each end of the \textit{$n_p th$} step are represented by the next higher approximation,

 \begin{multline}
 \delta \theta_{n_p} \equiv \delta \theta _1 = \pi - 4 \pi p +\frac{2 \pi}{n_p}(p - 1/2)^2  \\
  \delta \theta_{n_p-1} \equiv \delta \theta _2 = \pi - 4 \pi p +\frac{2 \pi}{n_p}(p + 1/2)^2
 \end{multline}
\noindent
 then the angles in Figure \ref{construct} are given by

 \begin{equation}
\theta _1 = \pi /2 +\delta \theta_1/2  \qquad
 \theta_2 = \pi/2 + \delta \theta_2/2  \qquad
 \theta_3 = \pi -\theta_1 -\theta_2
 \end{equation}

 In this case, application of the \textit{Sine Law} for adjacent sides and angles gives

\begin{equation}
P(s) = \sum_1^{n_p} n^{-s} + \frac{sin(\theta_2)}{sin(\theta_3) n_p^{1/2}}e^{i(-t log(n_p) - \pi /2 +\delta \theta_1 /2)}
\end{equation}

Pendant centers found by this construction are indicated in Figure \ref{pend3x3}.  Note the red $\mathcal{+}$ in the upper left of Figure \ref{pend3x3}, which is the ``center'' of the pendant in the lower left having $\delta \theta = 0$; even the higher order correction can lead to errors whenever the vertex angle in Figure \ref{construct} goes to zero.  Note, however, that the point located as the center lies on the axis of symmetry.  This is related to the fact that zeros in the denominator of the correction term are accompanied by zeros in the numerator in Riemann's result, to be discussed below.

One procedure to avoid divergences in the calculation of $L$ is simply to calculate the angle differences using the definitions of Equation 5. Another is to calculate the displacement $L$ for steps adjacent to $n_p$ and choose the smallest.  This avoids anomalously long distances due to the effects above.  Also, there is a prescription for the maximum possible $L$, based on the radius of the pendant circle as a function of $(n_p,p)$.  Using the DTS with $\delta \theta = \pi -4 \pi p$, $\delta^2 \theta = 2 \pi \to 0$ and $\delta^3 \theta = -4 \pi /n_p$ and solving for the number of steps giving an angle change of $-\pi$, leads to the cubic equation

\begin{equation}
\delta n^3-3 \delta n^2 + (2+6 n_p (1-4p))  \delta n = 6 n_p
\end{equation}

When $p=0$, this gives a $\pi$-reversal in one step, which is obvious since $\delta \theta = \pi$.  However, when $p=1/4$, $\delta \theta = 0$ and $\delta n \simeq (6 n_p)^{1/3}$.  This gives

\begin{equation}
L \simeq 6^{1/3} n_p^{1/6}/2\pi
\end{equation}

Clearly, this result only holds for extremely large $t$.  As will be shown below, the error in $L$ has a small impact on the determination of $\zeta$ due to a cancelation effect, but is important when investigating the distribution of $P(s)$.

Figure \ref{pend6710} shows the actual results for $t=7007.189$, for which $n_p = 33$ and $p=.3950$.\footnote{This value corresponds to \textit{Gram point 6710}. The \textit{nth} Gram point is the value of $t$ for which $t/2 \log{t/2\pi}-t/2-\pi/8 = n \pi$.}.  The vertex marked $\oplus$ is the end of step $n_p$ and the beginning of step $n_p +1$, with the steps progressing in the negative (clockwise) sense.  The points $P_n$ for steps $n_p -1, n_p, n_p+1, n_p+2$ are denoted by $o,\times,+,\diamond$. Note that $N_{n_p}$ and $N_{n_p+1}$ are nearly coincident, defining the center of an otherwise nearly unrecognizable pendant.

 \begin{figure}[pt]
\centering
\includegraphics[angle=0,width=.6\textwidth]{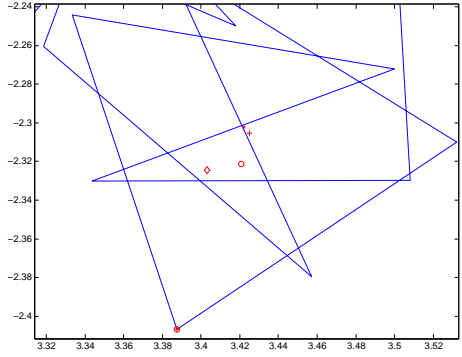}
\caption{Steps near the pendant center of Gram point 6710, with symbols defined in the text.}
\label{pend6710}
\end{figure}

\it
\item For Euler spirals with sufficiently large $M$, the sum over the steps between successive scroll centers is approximated by an integral of gaussian form whose value is

\begin{equation}
\frac{n_p^{1-2\sigma}}{k^{\sigma-1}} e^{i(-t log(t/2\pi) +t +(1-(k/n_p)^2)\pi/4}
 \end{equation}
 \rm

 The demonstration of conjugacy above, using a careful definition of the spiral centers, was purely numerical. When $M=(n_p/k)^2$ is sufficiently large, it is also possible to derive an expression for the distances between adjacent spiral centers of the \textit{kth} Euler spiral.  Application of the DTS leads to a complex integral of gaussian form:

   \begin{equation}
   \sum_L^U n^{-s} \simeq \frac{e^{i\theta_k}}{n_k^{\sigma}} \int_{-\infty}^{\infty} \frac{e^{i/2\delta^2\theta_k x(x-1)}}{(1+x/n_k)^{\sigma}} d x
   \end{equation}
\noindent
where the lower limit $L=[t/(2k+1)\pi]-1/2$ and the upper limit $U=[t/(2k-1)]-1/2$.  The notation ``-1/2'' signifies subtraction of 1/2 of the step, as in the zero order determination.  Although the details of the sum limits are stated for completeness, this are not important, as it will be assumed that they are $\pm \infty$.  The discrete variable $\delta n$ has been replaced by $x$ and it is assumed that it can be accurately represented by an integration. Expanding the denominator, the first order term is an odd function giving zero upon integration, and the second order term is $\mathcal{O} (n^{-(2+\sigma)})$\footnote{The Euler-Maclaurin integral approximation for the discrete sum is not appropriate unless $n > t/\pi $, where the straight discrete steps are sufficiently close to the arc of the smooth function $n^{-s}$ that a Bernoulli polynomials can correct for the error (see Figure 5.).  The DTS, however, is valid for this integration.  Error arises primarily when $M$ is not sufficiently large for the integral approximation to the sum.}.

The DTS expansion to $\delta^2 \theta$ is not applicable over the entire range $(L,U)$, but because the gaussian is narrow, it does not need to be.  In figure \ref{squaresteps}, the real and imaginary step sums are plotted vs. step number.

 \begin{figure}[pt]
 \centering
 \includegraphics[angle=0,width=.6\textwidth]{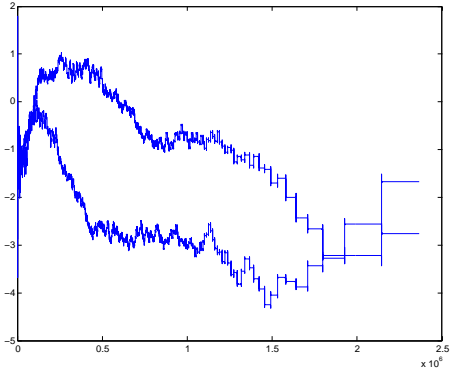}
\caption{Real and Imaginary Components of $\sum^N n^{-s}$ vs. N, $s=1/2+10^9i$}
\label{squaresteps}
\end{figure}

Note that the angle change is localized to steps near the inflection points of the Euler spiral.  Alternating directions of the small steps within the spiral centers rapidly average to constant values.\footnote{This linear step plot also shows the large fraction of the total number of steps in a few final Euler spirals.  Since each of these is equal to a single step under the substitution $s \to 1-s$, the numerical effectiveness  of the Riemann-Siegel sum, to be discussed below, is clearly shown. Since this integration primarily illustrates features which will be replaced by more precise analytic results, detailed error analysis is omitted.}

The center, or inflection point, of the \textit{kth} Euler spiral is defined as the integer step where $\delta \theta = 2\pi k$, or

\begin{equation}
 n_k = [(e^{2\pi k/t}-1)^{-1}]
\end{equation}

As in the first order spiral center location, define $\delta t \in [0,2\pi k)$ by $t-2\pi k n_k$.  As $\delta t$ varies from $0 \to 2\pi k$, the step angle at $n_k$ changes by $2 \pi k$.  The center angle at step N is then

 \begin{equation}
 \theta_k = -t log\frac {t}{2k\pi} -2k\pi N +t
 \end{equation}

The middle term is zero by the mod function.  At the kth Euler spiral, the second angle difference is $2\pi(k/n_p)^2$. Completion of the square, transformation to $y=\sqrt(-i)(k/n_p)(\delta x -1/2)$, and use of the value $\sqrt{2\pi}$ for the gaussian integral gives;

\begin{equation}
\frac{n_p^{1-2\sigma}}{k^{\sigma-1}} e^{i(-t log(t/2\pi) +t +\pi/4(1-(k/n_p)^2)}
 \end{equation}

The ``scroll tilt angle'' at the center is $\pi /4$. For ease of integration, the centered scroll angle was used instead of the forward step angle. The difference between the forward and centered angles is $1/2 (1/2)(1/2-1) \delta^2 \theta$ from the DTS, and therefore the term $-\pi/4 (k/n_p)^2$ is properly omitted when the forward difference angle is used.  For $k$ within the range of applicability, this scroll integral gives the relation between a conjugate region and the original step with $s \to 1-s$. As noted above, this factor with $k=1$ is the ratio $Q(s)$ defined above.

It is now possible to examine the geometric origins of the angle term $-2\Theta (t)$ in $Q(s)$\footnote{This definition agrees with the standard notation of analytic number theory.  The significant terms in \cite[Eqn.1,Sec6.5]{edwards} are all describable geometrically, the higher order terms arising from more precise expressions for the $\Pi$-function and the gaussian scroll integral.}. The term $-t log(n_p^2)$ is $n_p^{-2 i t}$, or twice the step angle at the center of symmetry.  The use of $\delta t = t - (2k-1)\pi N$ to correct for the discretization to integer step values gives the term $it$, since the \textit{reference} term in $\delta t$ is zero mod $2\pi$.  The additional term $\pi/4$ corrects for the inherent ``tilt angle'' of a conjugate region, seen from the final Euler spiral of Figure \ref{xy9} with $M \simeq 2 \cdot 10^8$ steps to the $M=5$ steps of Figure \ref{pair5}.  Note that the amplitude of $Q$ varies monotonically from $1/n_p \to n_p$ across the critical strip, and is only equal to 1 when $\sigma = 1/2)$.

Because the sum of the angle of each original step and that of its conjugate region is $2\Theta$, the angle of the axis of symmetry of the complete Argand diagram is given by $\Theta(t)/2 - \pi /4$, or

\begin{equation}
\theta_{sym}(t) = -\frac{t}{2} log \frac{t}{2\pi} + \frac{t}{2} -\frac{3 \pi}{8}+\ldots
\end{equation}

 \it

\item There is a natural ``end'' to the series at step $[t/\pi-1/2]$.

The reciprocity relation between original and conjugate steps, with the symmetry condition, implies a finite series extending from the origin to its conjugate step.  Thus, after step $N=[t/\pi - 1/2]$, a final spiral develops whose radius shrinks to zero at the point $\zeta (s)$ for $\sigma > 1$ or grows without bound for $\sigma < 1$.

\rm

Consider a pair of adjacent steps $n,n+1 \gg t/\pi$.  Their perpendicular bisectors intersect at a distance $R = n^{-\sigma}/ \delta \theta \simeq n^{1-\sigma }/t$, defining a radius of curvature of the spiral arc that decreases with $n$ for $\sigma > 1$.  In this case, the final spiral collapses to the value of $\zeta$.  When $\sigma < 1$, $R$ increases without bound, as is well known for the Dirichlet series with a pole at $s = 1$.

At this point, the essential numerical and geometric aspects of the geometry of $\zeta$ have been described, and detailed comparisons with the results of ``traditional'' analytic number theory can be discussed.
\\

\end{enumerate}

\section{Summary of Analytic Number Theory Results}

The enumeration of comparisons between geometric and analytic results is preceded here by a list of significant ``classical'' results.  This subject is much too extensive to be covered briefly, and only the principle features are described, selected to be directly applicable to this discussion.  A familiarity with the references is necessary for full understanding\cite{apostol,edwards}.

\begin{enumerate}
\it
\item{Euler-Maclaurin Summation}
\rm

The sum for $\zeta (s)$ can be evaluated up to some $N>t/2\pi$, with the remaining terms being approximated by \textit{Euler-Maclaurin summation}\cite[sec.~6.2]{edwards}.  Approximating the straight steps by a integral over the smooth function gives:

\begin{multline}
\zeta (s) = \sum_1^N n^{-s} -\frac{N^{-s}}{2} + \frac{N^{1-s}}{s-1} + \frac{B_2}{2} s N^{-(s+1)}
+ \ldots \\ \frac{B_{2k}}{(2k)!)}s(s+1) + \ldots (s+2k-2)N^{-(s+2k)-1} \ldots + R_{2k}
\end{multline}
\noindent
where $B_n$ is a \textit{Bernoulli number}.

\it
\item Riemann's Contour Integration
\rm

The zeta function can also be expressed as an integral in the complex plane of the integrand

\begin{equation}
 \frac{(-x)^s}{(e^x-1)}\frac{dx}{x}
\end{equation}

 Three contour paths used in the analysis are shown in Figure \ref{contour} \cite[Chap.7]{edwards}:

\begin{figure}[pt]
\centering
\includegraphics[angle=0,width=1.2\textwidth]{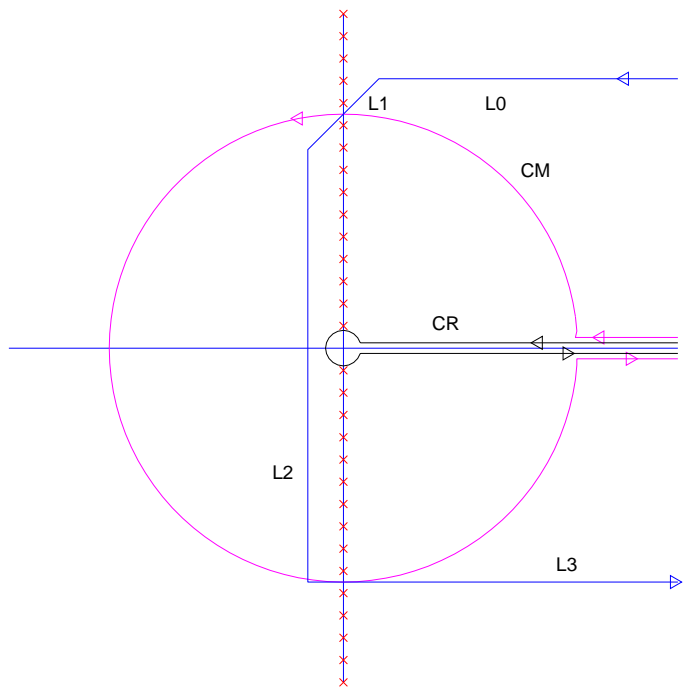}
\caption{Riemann's contour (CR), a contour also enclosing multiple poles on the imaginary axis (CM), and a segmented linear contour for saddle-point integration ($L_n$)}
\label{contour}
\end{figure}

Riemann's expression giving analytic continuation for the integrand $(-x)^s/(e^x-1)dx/x$ creates a ``cut'' on the positive real line by defining $-1$ in $(-x)^s$ as $e^{i\pi}$ on the upper horizontal leg and $e^{-i\pi}$ on the lower, return leg.  He argued that the circular loop of radius $\delta$ around the origin vanishes as $\delta \to 0$, leaving

\begin{equation}
\int_{CR}\frac{(-x)^s}{e^x-1}\frac{dx}{x} = lim_{\delta \to 0}\bigg \{ \int_{+\infty}^\delta \frac{(e^{i\pi} x)^s}{e^x -1} \frac{dx}{x} + \int_{\delta}^{\infty}\frac{(e^{-i\pi} x)^s}{e^x -1} \frac{dx}{x} \bigg\}
\end{equation}

The traditional procedure uses the expression $2 i sin(\pi s) =(e^{i \pi  s}-e^{-i \pi s})$, together with an identity for the factorial function, to define $\zeta$ directly in terms of the integral around CR:

 \begin{equation}
 \int_{CR}\frac{(-x)^s}{e^x -1}\frac{dx}{x}=2isin(\pi s) \Pi(s-1)\sum_{n=1}^{\infty}\frac{1}{n^s}
 \label{usualdef}
 \end{equation}

Although this was defined for $\sigma > 1$, it ``remains valid'' otherwise, ``since the integral clearly converges for all $s$, and the function it defines is complex analytic''.

 \it
\item Individual Steps
\rm

Since $1/(e^x -1)$ can be expanded as a sum of terms $ e^{-mx}$, it is of interest to present the result of integration for one of these over ``Riemann's contour'' CR

\begin{equation}
\int_{CR} (- x)^s e^{-mx} \frac{dx}{x} = (-1)^s \Pi(s-1) m^{-s}
\end{equation}

Using the Stirling approximation for $\Pi(x)$ provides an explicit expression in polar form for the coefficient of $m^{-s}$ in the integral

\begin{multline}
(-1)^s \Pi(s-1) \to \\
e^{-2\pi i (\sigma+it)} \sqrt(2\pi) (s-1)^{(s-1/2)} e^{1-\sigma -it} \to \\
\{ e^{\pi t/2} (2\pi)^{1/2} t^{\sigma-1/2}\}exp \{i[t log t -t - \frac{\pi}{4} -\frac{\pi \sigma}{2}] \}
\end{multline}

The $\zeta$ function is therefore expressed as the sum over all the individual terms $m^{-s}$, divided by the above coefficient

\begin{equation}
\zeta(s) = \frac{1}{(-1)^s \Pi(s-1)} \int_{CR} \frac{(-x)^s}{e^x -1} \frac{dx}{x}\
\end{equation}

\it
\item Integration of the Poles
\rm

The poles along the imaginary axis, indicated on Figure \ref{contour}, are integrated in conjugate pairs at $x=\pm 2 \pi i m$ (although only one pole of each pair is appreciable, depending on the sign of $t$).

\begin{equation}
(2\pi i)(-1)^s(2\pi i m)^{s-1}
\end{equation}
\noindent
which is expressed in polar form as

\begin{equation}
\{ e^{\pi t/2} (2\pi)^{\sigma} \} e^{i\{t log 2\pi -\frac{\pi \sigma}{2} + \pi \} }m^{s-1}
\end{equation}

Note that the integral around the pole pair has the same general form as an individual term in Riemann's contour, with $\Pi(s-1) m^{-s}$ replaced by $(2\pi i m)^{s-1}$.

\it
\item{Saddle-Point Integration}
\rm

The contour ``CM'' encloses the poles up to $\pm 2\pi i M$, as well as the discontinuity on the positive real axis.  When the first M terms of the sum are evaluated explicitly, the integrand is changed to $(-x)^s e^{-M}/(e^x-1) dx/x$ in order to evaluate the terms beyond M.  The derivative of the integrand is zero at the ``saddle point'' $x= (\sigma -1 + it)/M$, which is almost exactly on the imaginary axis at $it/M$ for large t \cite[sec. 3-6]{mathews&walker}.  The contour ``$L_n$'' is equivalent to ``CM'', since no poles or discontinuities of the the integrand are crossed as the contour is deformed from $CM \to L_n$.  However, this allows explicit integration along the contour, as the line ``L1'', crossing the saddle at an angle $\pi/4$ from upper right to lower left, is the only segment of ``$L_n$'' having a significant value of the integrand.    On that contour segment, the function approaches a narrow gaussian similar in form to the integral over an Euler spiral;  the angle of this segment is due to the transformation to the gaussian ``saddle''.  The result of this ``saddle point'' or ``steepest descent'' integration is $(-1)^s \Pi(s-1)$, which is the same coefficient as in the integration of $e^{-mx}$, here multiplying $N^{-s}/2$ instead of $m^{-s}$.

\it
\item  Landau Sum
\rm

The sum of $x^\rho$ over roots $\rho = \alpha_i + it$ has an interesting relationship to the primes\cite{landau}:

\begin{equation}
\sum_{pairs \rho _i} x^{\rho} = - \frac{T}{2\pi} \Lambda (x) + \mathcal{O}(log T)
\end{equation}
\noindent
where T is the largest $\alpha$ in the sum and \textit{Mangoldt's function} $\Lambda (x)$ is 0 unless x is a prime or a the power of a prime, where it equals $log p$.  This is the relationship shown in the introduction, except for the x-dependence of the numerator in Riemann's expression.  It will be discussed with respect to the average of real steps below.

\it
\item  Other Forms of the Zeta Function
\rm

Current mathematical research involves generalizations of Riemann's $\zeta$ function\cite[Chapter~12]{apostol}  The \textit{Hurwitz} $\zeta$ function  $\zeta(s,a)$ begins at the fixed real number $a$;

\begin{equation}
\zeta(s,a) = \sum_0^ \infty \frac{1}{(n+a)^s}
\end{equation}

When $a \ne 1$, the bilateral symmetry is no longer present.

Dirichlet L-functions are sums of Hurwitz $\zeta$ functions weighted by the characters of a finite abelian group. This is the \textit{the group of reduced residue classes modulo a fixed positive integer $k$.} The character table for $k=7$ is shown in Table \ref{characters}.  Dirichlet used this function in his theorem on primes in arithmetic progressions.

\begin{table}[h]
\centering
\begin{tabular}{|c|c|c|c|c|c|c|c|}
\hline
n & 1 & 2 & 3 & 4 & 5 & 6 & 7 \\
\hline
$\chi_1(n)$ & 1 & 1 & 1 & 1 & 1 & 1 & 0 \\
$\chi_2(n)$ & 1 & 1 &-1 & 1 &-1 &-1 & 0 \\
$\chi_3(n)$ & 1 & $\omega^2$ & $\omega$ & $-\omega$ &-$\omega^2$ &-1 & 0 \\
$\chi_4(n)$ & 1 &$\omega^2$ & -$\omega$ &-$\omega$ & $\omega^2$ & 1 & 0 \\
$\chi_5(n)$ & 1 &-$\omega $ & $\omega^2$& $\omega^2$ &-$\omega$ & 1 & 0 \\
$\chi_6(n)$ & 1 &-$\omega$ &-$\omega^2$ & $\omega^2$ & $\omega$ &-1 & 0 \\
\hline
\end{tabular}
\caption{Character table for k = 7, $\omega = e^{i \pi /3}$}
\label{characters}
\end{table}

\begin{equation}
L(s,\chi) = \sum_1^\infty \frac{\chi(n)}{n^s} \to \sum_1^k \chi (r) \zeta(s,\frac{r}{k})
\end{equation}

\end{enumerate}

\section{Connections of Geometric and Analytic Results}

It is now possible to relate the two views of $\zeta$.  First, note that symmetry implies a finite series from the first step to $[t/\pi]$.  The interpretation of \textit{geometrical continuation} as truncation of the infinite series at the point conjugate to the origin is introduced here.  It has effects on other comparisons, therefore it is appropriate to examine it first.  In particular, the finite series can be ``divided'' at the symmetry point into original steps and their conjugate regions, as discussed above and described below using complex analytic methods.

The infinite series representation of a function may increase without bound even in regions where the function from which it was derived is well-behaved.  This could be due to a convergence condition for the series that indicates divergence for \textit{all} radii around the expansion point once one position on a bounding circle has a pole.   The standard methods for analytic continuation of a series that would otherwise grow without bound are based on processes that represent the sum in ``other ways'', such as: solving the Cauchy-Riemann equations using boundary values in a region where the function is well-behaved; representing the series sum by an integral approximation and using partial integrations to represent a converging sum; ``walking'' a series of convergent regions around a singularity; or manipulation of conditionally convergent series by regrouping terms\cite[Sec.~1.4]{edwards}.  Riemann characterized these methods generally as ``finding a solution which is valid where the series converges, and which `remains valid' in other regions''.  In these sections, a case is made that appropriate truncation of the infinite series meets Riemann's standard, using three approaches:

\begin{enumerate}
\item

The geometric observation that the center of the final spiral has the desired properties.  While this is distinctly different from the above approaches, it will be supported by numerical examples and analytic comparisons.
\item

 The geometric determination of the spiral center agrees with the results of the Euler-Maclaurin formalism.  This is demonstrated to first order, and a method to incorporate higher order terms is indicated.
\item

Finally, it is shown that Riemann's contour integral leads to a limit process in which the incomplete gamma function attenuates step growth of the final spiral.  This occurs  in such a manner that the infinite sum approaches the spiral center.
\end{enumerate}

This procedure is admittedly quite different from the traditional approach; it is hoped that these detailed steps will convince the reader that the truncated series with appropriate end correction is a mathematically proper alternate prescription for continuation.

\begin{enumerate}
\it
\item Geometric Continuation using series truncation
\rm
   Figure \ref{circles} shows a portion of the ``final spiral'' beyond $[t/\pi]$.  This is a sequence of $32 \pi$ final steps ending at $n=549$ for $t=549.4975\ldots$, a zero of $\zeta$ when $\sigma = 1/2$. The values of $\sigma$ range in steps of $.05$ from $0.2$ to $1.4$. The zero at $\sigma = 1/2$ is marked by a red diamond. Values of $\zeta(\sigma +it)$, determined by explicit summation, are plotted in red for $\sigma \le 1$ and in green for $\sigma > 1$.  The spirals for $\sigma > 1$ are collapsing to the correct value where the series converges.  Those for $\sigma \le 1$ are expanding, but remaining centered on the continued value of $\zeta$.

 \begin{figure}[pt]
\centering
\includegraphics[angle=0,width=.6\textwidth]{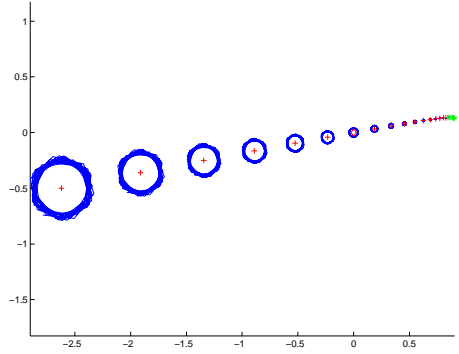}
\caption{Final spirals for $t=549.4975 \ldots$ and a sequence of values $\sigma \in [.2,1.4]$}
\label{circles}
\end{figure}

While this figure illustrates Riemann's concept of continuation, a quantitative argument uses the zero and first order locations of the final spiral center $N=[t/\pi]-N^{-s}/2$ evaluated for zeros of $\zeta$.  In this case, the value \textit{is} the error.  Figure \ref{error0} shows the zero order final spiral center $\sum^N n^{-s}-N^{-s}/2$ for $N=[\alpha_i/\pi]$, with small flags showing the direction of the final step for each of the 82 values of $\alpha_i \in (1100,1200)$.  The transverse and (small) longitudinal errors are $\mathcal{O}(N^{-3/2})$.

\begin{figure}[pt]
\centering
\includegraphics[angle=0,width=.6\textwidth]{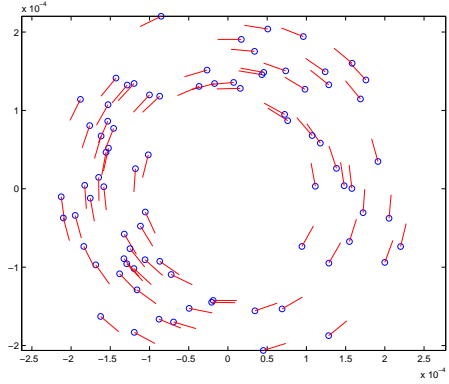}
\caption{Zero order errors for 82 zeros of $\zeta$ with $\sigma = 1/2$ and $\alpha_i \in (1100,1200)$}
\label{error0}
\end{figure}

The same calculation using first order corrections gives Figure \ref{error1}, in which the errors are reduced to $\mathcal{O}(N^{-5/2})$, the next higher order of omitted terms.

\begin{figure}[pt]
\centering
\includegraphics[angle=0,width=.6\textwidth]{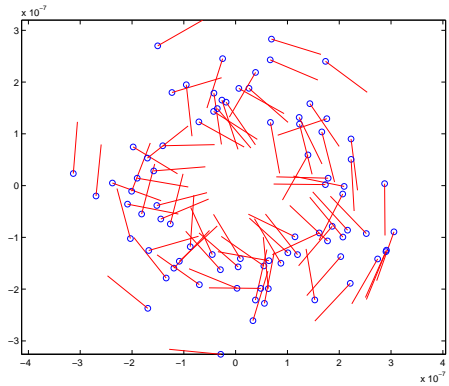}
\caption{First order errors for 82 zeros of $\zeta$ with $\sigma = 1/2$ and $\alpha_i \in (1100,1200)$}
\label{error1}
\end{figure}
\it
\item Geometric Continuation using Euler-Maclaurin summation
\rm

The Euler-Maclaurin expression given above can be rearranged into a form more appropriate for this application:

\begin{multline}
\zeta (s) = \sum_1^N n^{-s} -\frac{N^{-s}}{2} + \frac{N^{1-s}}{s-1} \sum_{k=0} \frac{B_{2k} s^{2k}}{(2k)! N^{-2k}} \\
+ \frac{N^{1-s}}{s(s-1)} \sum_{k=1} \frac{k(2k-3)B_{2k} s^{2k}}{(2k)! N^{-2k}} + \ldots
\end{multline}

 The zero order correction has been explicitly separated, leaving the higher order corrections to the sums over $k$. When the sum is carried to a sufficiently large value of $s/N$, it can be truncated at the first term,

 \begin{equation}
  \zeta (s) = \sum_1^N n^{-s} -\frac{N^{-s}}{2} + \frac{N^{1-s}}{s-1}
  \end{equation}

  In this case however, the sum must be evaluated at $t/\pi$, which is close to the convergence limit $t/2\pi$.  In addition, there is the ``binomial product'' $(s-1)s(s+1)\ldots$ giving decreasing powers of $s$. The monic polynomial resulting from multiplication of this product shows that the first term or two suffice:

\begin{equation}
(s-1)s(s+1) \ldots (s+2k-2) n\\ \to s^{2k} + k(2k-3)s^{2k-1} + \ldots
\end{equation}

\noindent when the ratio $2k^2/t \ll 1$.

The sum over powers of $(s/N)$ converges slowly for $\lvert s/N \rvert >\simeq \pi$, but it can be represented by the Bernoulli number generating function.   In this way, the first sum becomes

\begin{equation}
\frac{N^{-s}}{2} \frac{s}{s-1} \frac{e^{s/2N}+e^{-s/2N}}{e^{s/2N}-e^{-s/2N}}
\end{equation}

The ratio $s/(s-1) \to 1$ for large $s$, and manipulation of the definitions of trigonometric and hyperbolic functions lead to

\begin{equation}
\frac{N^{-s}}{2} \frac{tanh(\sigma/2N)-i\cdot cot(t/2N)}{1-i \cdot tanh(\sigma/2N) cot(t/2N)}
\end{equation}

 The first order expansion of the real term gives the longitudinal correction term $\lambda_N$. The term $cot(t/2N)$ implements the mod function, and its expansion about the zero of the $cot$ gives $\tau_N$, resulting in the same expression as the geometric first order correction.

Expansion of the denominator of the first sum leads to corrections of $\mathcal{O} (N^{-3-\sigma})$. Therefore, it is the \textit{second} sum term (together with second order terms from the first sum) that contains the second order correction; here it is evaluated at $t/N = \pi$ to determine an estimate of the first order error.  This leads to the interesting expression

\begin{equation}
\frac{N^{-s}}{2} \sum_{k=1} \frac{k(2k-3)}{2^{2k}} \zeta (2k)
\end{equation}

using Euler's famous relation\cite[Sec.~1.5]{edwards}

\begin{equation}
\zeta(2k) = \frac{(2\pi)^{2k}(-1)^{k-1}}{(2k)!} B_{2k}.
\end{equation}

Because this sum converges rapidly (the ratio test gives $term_{k+1}/term_k \rightarrow 1/4$ using $\lim {\rightarrow \infty}\quad \zeta(2k) = 1)$, the constant implicit in $\mathcal{O}(N^{-2-\sigma})$ is well-behaved.

\it
\item Geometric Continuation using Riemann's Contour Integral
\rm

This is addressed by using $N= [t/\pi]$ in the saddle-point integration.  The lowest order of the gaussian integral gives the zero order correction $-N^{-s}/2$.  Higher order terms could be derived by following Riemann's method for extension of the accuracy in the Riemann-Siegel equation, discussed below.  Since this must lead to the corrections already sketched for the Euler-Maclaurin method, it will be left for further examination by the reader.  A far more fundamental observation relates to the limit as the radius $\delta$ about the origin in Riemann's contour approaches 0.

  Referring to Eqn. \eqref{usualdef}, note that, depending on the sign of $t$, the factors $e^{i\pi s}$ and $e^{-i\pi s}$ are remarkably different in magnitude.  Their ratio is $\simeq e^{2\pi t}$, and even before the limit is taken, only the larger term remains.  This yields

 \begin{equation}
  \int_{CR}\frac{(-x)^s}{e^x -1}\frac{dx}{x}= lim_{\delta \to 0} e^{-i\pi s} \int_{\delta}^{\infty}\frac{x^s}{e^x -1} \frac{dx}{x}
  \label{newdef}
  \end{equation}

Combining the two representations of the integral around CR:

\begin{equation}
\Pi(s-1) \zeta(s) = lim_{\delta \to 0}\int_{\delta}^{\infty}\frac{x^s}{e^x -1}\frac{dx}{x}
\label{dampedspiral}
\end{equation}

Retaining the nonzero lower limit, the variable transformation $y = mx$ gives the previous result with an additional factor $\Phi(s-1,m\delta)$ in the coefficient\footnote{The notation Q is used in \cite{recipes} for this factor, where it is noted that it is unique to that reference. This work uses $\Phi$, since Q is used here for the complex factor relating initial and conjugate steps.  For consistency, the expression ``incomplete $\Pi$ function'' would be proper, but will not be used here.  The context provides any necessary changes to arguments used in the references.},

\begin{equation}
\zeta (s) = \sum_1^{\infty}\frac {\Phi (s-1,m\delta)}{n^s}
\end{equation}

where $\Phi(s-1,m\delta)$ is the \textit{regularized upper incomplete gamma function}

\begin{equation}
\Phi(s-1,m\delta) = \frac{\Pi(s-1,m\delta)}{\Pi(s-1)}
\end{equation}

\noindent and

$\Pi(s-1,m\delta)$ is the \textit{upper incomplete gamma function}\cite[Sec.~6.2]{recipes}

\begin{equation}
\int_{m\delta}^\infty x^{s-1} e^{-x} dx
\end{equation}

$\Phi$ approaches 1 for $m\delta \ll 1$, and goes to 0 as $m\delta$ increases to $\infty$.  Since the peak of the integrand is at $x=s-1$, $m\delta \simeq t$ gives $\Phi \simeq 1/2$.  Calculated values for real arguments are plotted in \cite{recipes}.

If $\delta$ is chosen for the explicit summation such that $\Phi(s-1,N\delta) \cong 1$, the terms beyond N can be approximated by the Euler-Maclaurin approximation as an integral

\begin{equation}
\sum_N^\infty \frac {\Phi(s-1,n\delta)}{n^s}   \to \int_N^\infty \frac{1}{x^s}\{ \int_{x \delta}^{\infty} y^{s-1} e^{-y} dy \}dx
\end{equation}

There is no closed algebraic expression for the integral representing the incomplete gamma function, but the sum of the remaining steps can be evaluated using integration by parts:

\begin{multline}
\frac{x^{1-s}}{1-s} \Phi(s-1,x \delta ) \arrowvert_N^{\infty} + \int_N^{\infty} \frac{x^{1-s}}{1-s} (x \delta )^{1-s} e^{-x \delta} dx \\
= \frac{N^{1-s}}{s-1} + \delta^{1-s}\int_N^{\infty}\frac{x^{1-s}}{s-1}  e^{-x \delta} dx
\end{multline}

In the limit $\delta \to 0$, consideration of the attenuation of steps beyond $N = t/ \pi$ due to the regularized upper incomplete gamma function extends the region of convergence to $\sigma <1$.  This interpretation of analytic continuation leads to the same value of $\zeta(s)$ as the principle of geometric continuation.

\it
\item Detailed Symmetry
\rm

In the introductory presentation of analytic results, it was noted that the integral over the \textit{mth} pole pair gives a factor $2\pi i (-2\pi im)^s /(-2\pi im)$ rather than $ \Pi(s-1)$.  The ratio of the coefficients of the powers of $m$ in the two integrals is therefore

\begin{equation}
-\frac{(2\pi i)^s}{\Pi(s-1)} = \frac{(2\pi i)^{1/2} (2\pi i)^{s-1/2}}{(2\pi)^{1/2}(s-1)^{s-1/2} e^{1-s}}
\end{equation}

Manipulation of these terms gives

\begin{multline}
-(\frac{it+\sigma -1}{2 \pi i})^{1/2-s} i^{1/2} e^{1-\sigma -it} \\ =(\frac{t}{2\pi})^{1/2-s}\{(1+\frac{1}{x})^x \}^{\frac{s-1/2}{x}}e^{s-1+i\pi /4}
\end{multline}

with $x=it/(\sigma -1)$.  Since $\lvert x \rvert >> 1$, the limit definition of $e \to (1+1/x)^x$ for large $x$ gives\footnote{A more precise representation of the factorial function than the Stirling approximation would give higher order terms in the imaginary exponent.}

\begin{equation}
Q(s)=\frac{(2\pi i)^s}{\Pi(s-1)} \to n_p^{1-2s} e^{i(t+\pi /4)}
\end{equation}

\it
The \textit{mth} pole pair integration gives $Q(s) m^{s-1}$, the extent of the region conjugate to an initial step $m^{-s}$. Thus, integration along the positive real axis, supplemented by integrating the pole pairs along the imaginary axis, provides the definitive analytical justification for the symmetry.
\rm

\it
\item The Functional Equation
\rm

The difference between each original step and its corresponding conjugate region, summed over all steps, goes to zero, so

\begin{equation}
\zeta(s) = Q(s) \zeta(1-s)
\end{equation}

This conclusion is reached in classical analysis by integrating around a contour which gives 0 as the difference between Riemann's contour and the sum over all pole pairs\cite[Secs.~1.6,~1.7]{edwards}.  The functional equation is thus a global expression of the detailed symmetry.  One important consequence of the functional equation is that the occurrence of a zero at $s=1/2 + \epsilon +it$ implies that there must also be one at $s=1/2-\epsilon +it$.

\it
\item The Riemann-Siegel Equation
\rm

Riemann used symmetry to calculate $\zeta$ more efficiently than an explicit sum of steps beyond $t/2\pi$\cite[Chap.~7]{edwards}. He used contour $L_n$ with $N= n_p$, expressing the corrections to $\sum_1^{n_p}n^{-s}$ and its conjugate in terms of a saddle point integral.  He chose not to publish the result however, and mathematicians doubted that he done it until Siegel examined his handwritten notes and publicized the achievement.  Riemann had not only determined the first order correction, but derived a hierarchy of coupled equations giving arbitrary precision. This story and his complete mathematical development are described in \cite[Sec. 7.1]{edwards}.

The geometry of his method is illustrated in Figure \ref{rslines} for $s=1/2+2220000.15i$.\footnote{Although the procedure is valid for any argument, this choice of $s$ reduces confusing overlaps and displays the features clearly.}  The green vector connects the origin to the pendant center $P(s)$.  This, plus its conjugate $P(1-s)$, in black, ends at the center of the final spiral at $\zeta (s)$.  The bisector of these two lines is the symmetry axis $\theta_{sym}$ at the angle

\begin{equation}
\theta_{sym} = -\frac{t}{2} log \frac{t}{2 \pi} +\frac{t}{2} - \frac{3}{8} \pi + \ldots
\end{equation}
\noindent
and the line perpendicular to this through the origin is at the angle $\Theta$.\footnote{In \cite{edwards}, this is referred to as $-\theta$.}

\begin{equation}
\Theta = -\frac{t}{2} log (\frac{t}{2 \pi}) +\frac{t}{2} +\frac{\pi}{8} + \ldots
\end{equation}

\begin{figure}[pt]
\centering
\includegraphics[angle=0,width=.6\textwidth]{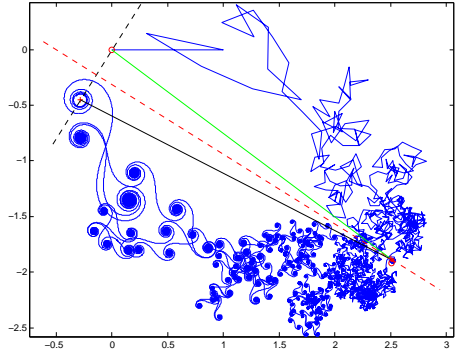}
\caption{Riemann-Siegel illustration for $s=1/2+2220000.15i$.}
\label{rslines}
\end{figure}

The significance of the Riemann-Siegel equation transcends efficient evaluation of $\zeta(s)$; it also leads to important geometric interpretations.\footnote{Discussions of the Riemann-Siegel procedure are often limited to $\sigma = 1/2$.  Because the angles are determined entirely by $t$, they are unchanged as $\sigma$ is varied from $1/2$.  The amplitudes $\mathsf{P}(s)$ and $\mathsf{P}(1-s)=n_p^{1-2 \sigma}\mathsf{P}(s)$ are to be used when $\sigma \ne 1/2$.}  Given $s = 1/2 + it$, the line through the origin at angle $\Theta$ can be constructed, and it will be noticed that  $\zeta(1/2+t)$ always lies on this line.  When $\sigma = 1/2$, $P(s)$ and $P(1-s)$ are of the same lengths and as t increases they superimpose to go through a zero. This changes the sign by $\pi$.  Therefore, its phase will be $\Theta$ or $\Theta - \pi$, depending on the number of zeros $N(t)$ encountered up to that value of $t$.

Another Argand representation of the geometry is shown for the Gram point used to illustrate the calculation of the pendant center, depicted in Figure \ref{gram6710all}.  The symmetry axis is drawn through $P(s)$ as determined previously (in green) and the $\Theta$-axis through $P(s)$ (in black).

\begin{figure}[pt]
\centering
\includegraphics[angle=0,width = .5\textwidth]{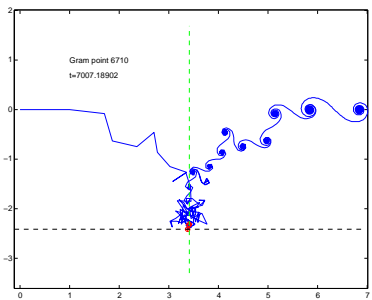}
\caption{Argand diagram for Gram point 6710, t = 7007.18902, with symmetry axes, $n_p$ and pendant center.}
\label{gram6710all}
\end{figure}

The enlargement of the pendant region in Figure \ref{6710axes} shows the geometry of the first order correction more clearly.

\begin{figure}[pt]
\centering
\includegraphics[width=.6\textwidth]{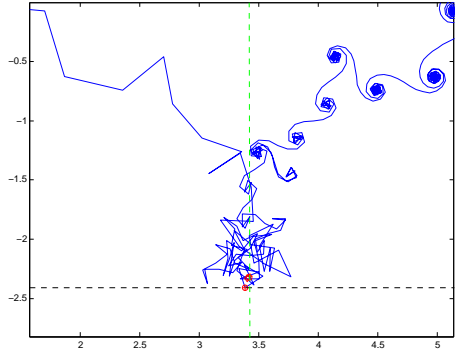}
\caption{Enlarged Argand diagram for Gram point 6710.}
\label{6710axes}
\end{figure}

The addition of the polar representations of $P(s) = \mathsf{P} e^{i\phi}$ and its conjugate yields

\begin{multline}
\zeta (s) = P(s) + Q P(1-s) \\ = \mathsf{P} e^{i\phi} + \mathsf{P} e^{2i\Theta - i\phi} \\
\to 2 cos(\phi -\theta) \mathsf{P} e^{i\Theta}
\label{rsequation}
\end{multline}

This expression shows once again that the complex phase of $\zeta$ is $\Theta$ or $\Theta -\pi$, depending on the sign of the cosine.

The geometric construction of the Riemann-Siegel method uses the pendant center $P(s)$ and the angle $\Theta$.  It is possible to make contact with the classical analysis, which used $\sum_1^{n_p}n^{-s}$ and a remainder term $R$.

 As discussed below in the section ``Zeros of $\zeta$'', the value of $\zeta$ is independent of the translation of $P(s)$ to \textit{any} point $P'$ lying on the symmetry axis.  In the classical development, $\sum_1^{n_p} n^{-s}$ and its conjugate $\sum_1^{n_p} n^{s-1}$ are both projected onto the $\Theta$ axis and their sum is termed $Z$, a real number.  This was done to simplify locating zeros by the sign change in this real quantity. The first order ``mismatch'' between $\sum_1^{n_p} n^{-s}$ and its conjugate was determined analytically by Riemann, as described above.  Geometrically, the vector from step $n_p$ (indicated by $\oplus$ in Figure \ref{6710axes} to $P(s)$ (marked $\mathsf{x}$) is $L = \mathsf{L} e^{i\theta_L}$.  If $L$ is also projected to the $\Theta$-axis, an expression for $\zeta$ results in terms of $\sum_1^{n_p} n^{-s}+R$, where $R = 2\mathsf{L} cos(\theta_L - \Theta)$.  The factor $2$ accounts for the contribution from each ``side'' of the sum of steps and conjugate steps to $n_p$.  Note that whereas the step angle in $\theta_L$ involves $t log(n_p)$, the $\Theta$-axis involves $t log(\sqrt(t/2\pi))=t log(n_p)+t log(1+p/n_p)$.  Substituting $2 \pi (n_p + p)^2$ for $t$, expanding the logarithm and discarding terms of $\mathcal{O}(1/n_p)$ gives

\begin{multline}
\Theta = -t log(n_p) +\pi n_p^2 + \pi (\frac{\pi}{8}- 2 p^2) \\
\theta_L = -t log(n_p) - \frac{\pi}{2} -\frac{1}{2} \delta \theta \\
\end{multline}

Noting that $\pi n_p^2 (mod 2\pi) = 0 $ for $n_p$ even or $ \pi$ for $n_p$ odd;

\begin{equation}
\theta_L -\Theta = 2 \pi (p^2 - p -1/16 + 0 \quad( or \quad \pi)  )
\end{equation}

Combined with the definition of $L$ above, this gives the pleasing result

\begin{equation}
R = (-1)^{n_p} \frac{cos 2 \pi (p^2 - p - 1/16)}{n_p^{1/2} cos2 \pi p}
\end{equation}

``Pleasing'', because this geometric construction agrees with the first order correction derived by Riemann from the saddle-point integration.  As described in \cite[Sec.~7.4]{edwards}, there are zeros in both the numerator and denominator of $R$ at $p=1/4,3/4$.  This can also be understood geometrically. The error in the pendant center location grows large at these values of $p$, but does so \textit{along} the symmetry axis, since when $\delta \theta = 0$ the $n_p$ step angle is parallel (or antiparallel) to the $\Theta$ axis.  It is therefore perpendicular to $\theta_L$.  When projected to the $\Theta$-axis, the zeros of $cos(\Theta - \theta_L)$ cancel this effect.

To demonstrate the first-order error numerically, the magnitudes of $\zeta (1/2 + i\alpha_k )$ were computed to lowest order by omitting the remainder $R$.  Results for values of $\alpha_k$ having $ k \in [97563,99596]$, all of which have $n_p = 108$, are plotted versus $p$ in Fig. \ref{rserror1}.  As in previous examples, the magnitudes of $\zeta$ at these zeros \textit{are} the errors, which agree with Riemann's calculation within a few times $10^{-4}$.  The next omitted correction term is $\mathcal{O}(1/n_p^{3/2}) \simeq 10^{-3}$.

\begin{figure}[pt]
\centering
\includegraphics[angle=0,width=.6\textwidth]{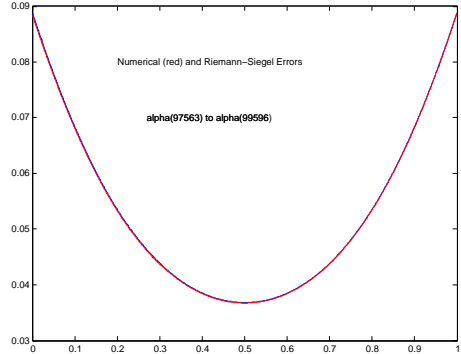}
\caption{Errors vs. $p$ for the lowest order calculation of $ \lvert \zeta \lvert $ for 2000 zeros of $\zeta$ having $n_p = 108$.}
\label{rserror1}
\end{figure}

Note that there is \textit{no} step node at which the first order error is zero.  When $\delta \theta \simeq 0$, the pendant center is circled by the Argand diagram steps, and when $\delta \theta \simeq \pi$ it is near the step midpoints.

\it
\item The Lima\c con
\rm

Plots of the Riemann-Siegel equation in the complex plane show a strong similarity to a classical geometrical figure called the \textit{lima\c con}.  Attributed to the father of Blaise Pascal, it is generated by the rotation of one line segment attached to the origin, with another segment pivoting about the end of the first.  In the usual applications, the lengths of the two segments are constant, but perhaps not equal, and the rotation rates are fixed.  In the Riemann-Siegel equation, the segment lengths are $\textsf{P(s)}$ and $e^{1-2\sigma}\textsf{P(1-s)}$, and the two angles are $\phi$ and $2\Theta - \phi$.  The lengths are thus equal when $\sigma = 1/2$, and zeros result when $\phi = 2\Theta - \phi (mod2\pi)$, as noted above.  In Figure \ref{limacon1000}, line segments from the origin to $P(s)$ (the dashed red curve), plus segments $QP(1-s)$, give $\zeta(s)$ (the solid blue curve)for Gram points $1000 \to 1004$.  Careful inspection of this figure provides a useful visualization of the Riemann-Siegel equation.

\begin{figure}[pt]
\centering
\includegraphics[angle=0,width=.6\textwidth]{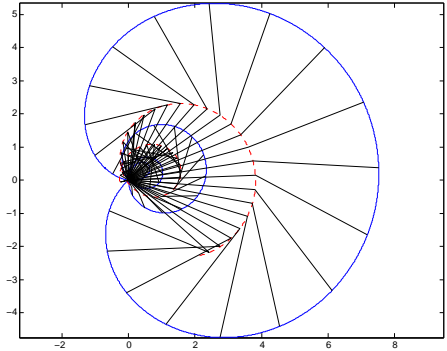}
\caption{Lima\c con construction for Gram points $1000 \to 1004$.}
\label{limacon1000}
\end{figure}

\it
\item Properties of P(s)
\rm

Figure \ref{psurf} shows the real and imaginary surfaces of $P(s)$ on a region of the critical strip from $t = 61425$ to $61429$.

\begin{figure}[pt]
\centering
\includegraphics[angle=0,width=.6\textwidth]{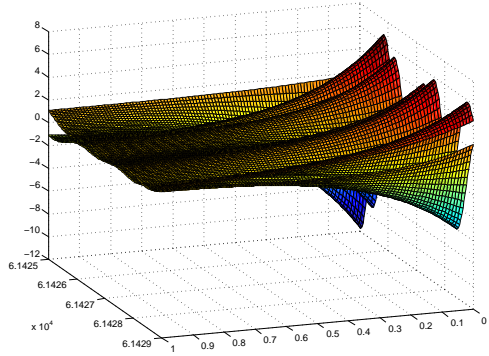}
\caption{$P(s)$ for $t \in [61425,61429]$, $\sigma \in $ [0,1]}
\label{psurf}
\end{figure}

For large $\sigma$, $P(s)$ approaches a real surface at 1 and an imaginary one at 0, just as for $\zeta(s)$. In fact, the general shape of $P(s)$ is nearly equal to $\zeta (s)$ when $\sigma$ is large, since they share the same steps up to $n_p$, and the step lengths decrease as $n^{-\sigma}$.

 Note the general tendency for the amplitude of the real and imaginary "flutes" to increase as $\sigma$ decreases.  Zeros of $P(s)$ occur at the intersections of the zeros of these real and imaginary "fluted surfaces".  This is illustrated in Figure \ref{zeroplots}, which shows zeros of $\mathcal{R} (P)$ in blue and $\mathcal{I} (P)$ in red.

 Also shown are zeros of $\mathcal{R} (\zeta)$ in black and $\mathcal{I} (\zeta)$ in green.  This figure shows that the ``symmetry zeros'' of $\zeta$ fall on the line $\sigma = 1/2$ as determined by Equation \ref{rsequation}. Zeros of $P(s)$ are seen in the figure as intersections of the blue and red curves.  For both $\zeta$ and $P$, the imaginary surface is ``centered'' on zero.  Therefore, the zeros of their imaginary components may extend across the entire critical strip.  The real components, however, are centered around 1, and lines of intersection with zero only give the parabola-like curves for smaller values of $\sigma$, where the ``real flutes'' have increased in amplitude.

 The functional equation arises when steps are paired with conjugate steps.  Because  $P(s)$ only includes steps up to the symmetry center, it is not subject to this symmetry.  Therefore, there is no condition for a zero at $s$ to imply a zero at $1-s$.  Figure \ref{zerop}, the Argand diagram for the zero of $P(s)$  shown in Figure \ref{zeroplots} at $s=.43+10010.8$, illustrates this situation.  Because $P(s)=0$, the pendant center and the origin are coincident, as indicated by the red $o$.  However, as indicated by the red $+$, $\zeta(.43+10010.8 i)$ is not zero.  In this case, $\zeta$ equals the second term in the Riemann-Siegel equation $n_p^{0.14}e^{2 i \Theta} P(1-s)$.

\begin{figure}[pt]
\centering
\includegraphics[angle=0,width=.6\textwidth]{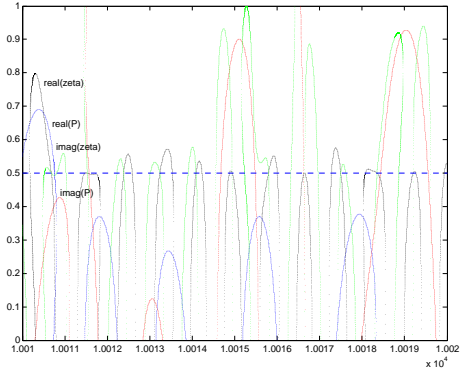}
\caption{Zeros of real and imaginary surfaces of $P(s)$ and $\zeta(s)$ on the critical strip for $t \in [10010,10020]$}
\label{zeroplots}
\end{figure}

 \begin{figure}[pt]
\centering
\includegraphics[angle=0,width=.6\textwidth]{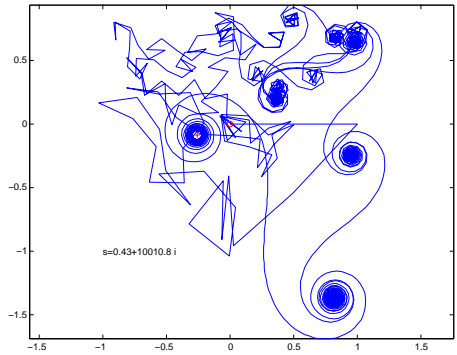}
\caption{Argand diagram for $s=.43+10010.8 i$, a zero of $P(s)$.}
\label{zerop}
\end{figure}

\it
\item Zeros of $\zeta$.
\rm

The popular literature concerning Riemann's $\zeta{s}$ is focused on the zeros in the critical strip and his surmise that they all lie on the line $\sigma = 1/2$.  The Riemann-Siegel equation facilitates the analysis of zeros, using the polar forms of $P(s)$ and $Q(s)P(1-s)$.  Their amplitudes are equal \textit{whenever} $\sigma = 1/2$, where $P(1-s)$ is the complex conjugate of $P(s)$ and $n_p^{1-2\sigma}$ equals one.  Although the $\sigma$ dependence of $P(1-s)$ is ``reversed'' across the critical strip, the factor $n_p^{1-2\sigma}$ counters this as shown in Fig. \ref{psurf120to129}.  For the single flute on the critical strip with $t \in (120,129)$, the two surfaces intersect transversally.

 \begin{figure}[pt]
\centering
\includegraphics[angle=0,width=.6\textwidth]{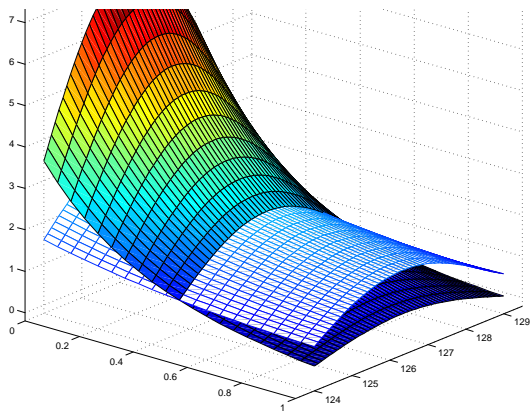}
\caption{Magnitudes of $P(s)$ (plain mesh) and $Q(s)P(1-s)$ (color) on the critical strip for $124<t<129$.}
\label{psurf120to129}
\end{figure}

 When the amplitudes are equal, a zero occurs when the angles satisfy

\begin{equation}
\phi = 2\Theta - \phi \to \Theta - \phi = (2k-1)\frac{\pi}{2}
\end{equation}

  Therefore, in addition to the distinction between zeros on the negative real axis and the \textit{nontrivial zeros} on the critical strip, the latter can be further categorized as either \textit{symmetry-zeros}, which are ``forced'' by the condition of symmetry, or those requiring equality (including zeros) for both $P(1/2- \epsilon +it)$ and $n_p^{2\epsilon}P(1/2+ \epsilon +it)$ for $\epsilon \ne 0$.

\it
 The symmetry zeros exist because the symmetry requirement forces a zero whenever $\sigma = 1/2$ and $\Theta - \phi$ is an odd multiple of $\pi /2$.

$\Theta$ is smoothly and monotonically decreasing function of $t$, but $\phi$ exhibits variations around its average trend, as shown in Figure \ref{phiphase}.  These variations account for the irregular spacing of symmetry zeros.
\rm

\begin{figure}[pt]
\centering
\includegraphics[angle=0,width=.6\textwidth]{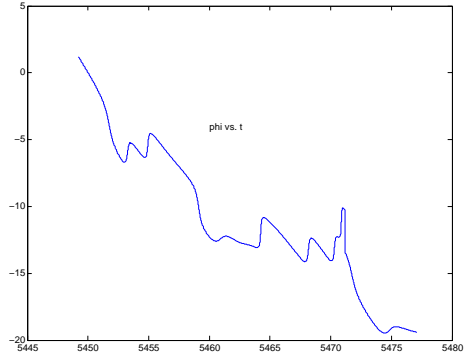}
\caption{Phase of $P(s)$ vs. t, $\sigma = 1/2$.}
\label{phiphase}
\end{figure}

It is not particularly easy to see the effect of $\phi$-variations from this figure.  When $mod (\phi -\Theta), 2\pi$ is plotted, the zeros, marked by $\mathsf{o}$, are seen to occur at odd multiples of $\pi/2$.  Figure \ref{phiminthet} shows the range of Gram points $n = 5000 \to 5008$ where $\Theta_n$ equals $n\pi$.

\begin{figure}[pt]
\centering
\includegraphics[angle=0,width=.6\textwidth]{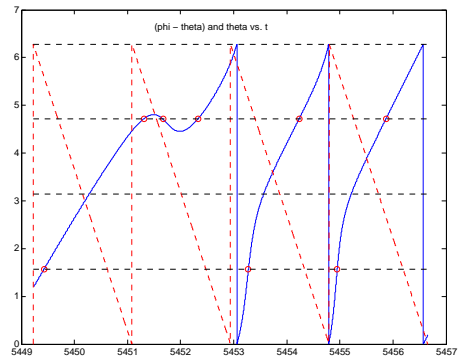}
\caption{$\Theta$ (red dashed line) and $(\phi - \Theta)$ vs. t, $\sigma = 1/2$.}
\label{phiminthet}
\end{figure}
\rm

The value of $\zeta$ is invariant under the translation $P(s) \to P'(s)$ along the symmetry axis when $\sigma = 1/2$.  This gives $\zeta = 0$ \textit{whenever} $P(s)$ lies on the line through the origin at angle $\Theta$.  This is the usual condition for a zero; the smooth advance of $\Theta(t)$ and the less regular variation of $\phi$ combine to give irregularly spaced zeros for any amplitude of  $\mathsf{P} \ne 0$ as $\phi \to \theta_{sym}$.  There does not appear to be any fundamental restriction on how closely $P(1/2+it)$ itself can approach zero, which would, however, still be consistent with Riemann's hypothesis. The general shape of the curves in Figure \ref{psurf} and the zeros in Figure \ref{zeroplots} indicate the challenge of determining simultaneous zeros of $P$ at $\sigma = 1/2 \pm \epsilon$.  Since no exceptions of RH have been reported to date, every zero of $\zeta$ that has been evaluated is a symmetry zero.

\it
\item Bounds
\rm

It is often of interest to establish bounds on the real and imaginary parts of $\zeta$. For example, what is the largest value of $\sigma$ for which $\mathcal{R}(\zeta)$ can be equal to zero? For $\sigma > 1$, the total length of the steps converges to the finite value $\zeta (\sigma + i0)$.  For example, Euler's equation gives $\zeta(4) = \pi^4/90 \to 1.082$.  Because the first step of the Riemann $\zeta$ function is always 1, it is easy to state that the largest possible magnitude of $\mathcal{I}\zeta$ is $i (\lvert\zeta(\sigma)-1 \lvert)$ and the smallest value of $\mathcal{R}\zeta(\sigma)-1$.  For $\zeta=0$, it must be that $\zeta(\sigma) \ge 2$.  This occurs near $\sigma = 1.7$.  These particular examples are quite elementary, but the approach may lead to new results.

Although the series of step lengths converges rapidly as $\sigma$ increases beyond 1, it is sometimes useful to have an upper bound for the magnitude of $\zeta$ within the critical strip.  Use of the Riemann-Siegel equation places limits on $P(s)$, and hence $\zeta$.  For example, $P(s)$ in Figure \ref{gram6710all} is about 3.5, and its maximum possible value is

\begin{equation}
\sum_1^{33} \frac{1}{n^{1/2}} \to 10.1
\end{equation}

The Argand diagram shows the predominant positive direction of the initial steps, which accounts for this uncharacteristically large value of $\zeta$, as seen in Figure \ref{6704}.

 It was stated qualitatively above that $P(s) \to \zeta (s)$ for large $\sigma$.  This can also be seen by the application of bounds;  for $\sigma = 0$, $P(s)$ is strictly bounded by the sum of $n_p$ unit steps.  (In practice, it is easily seen that this is far too generous, as $P(0+it)$ rarely exceeds ``a few''.)\footnote{Furthermore, the strict bound on $P(1+it)$ is $\sum_1^{n_p} 1/n \simeq log(n_p)$.}  Applying these to the Riemann-Siegel equation for $\zeta$;

\begin{multline}
|\zeta(1+it)| = |P(1+it) + np^{-1}e^{2i\Theta}P(0+it)| \\  \to |P(1+it) +e^{2i\Theta}|
\end{multline}

\it
\item Average values

The sum of an ensemble of Argand diagrams  $\sum n^{-s}$ over the $\Delta N \simeq ( \Delta T/2\pi log T/2\pi$  zeros $\alpha_i$ of $\zeta$ on an interval $\Delta T \ll T$ reveals that the real components only change appreciably at prime steps. At step $p$, the changes average to

\begin{equation}
-T log p/(2 \pi p^{1/2})
\end{equation}

The special nature of this average must be stressed; it is only the average over Argand diagrams giving $\zeta = 0$.

It is often noted that the real component of $\sum_1^n n^{-s}$ is generally larger than the imaginary component, because the first step is always $(1+i \cdot 0)$. The real part is 1 at the first step, and goes to zero at the final step for $t = \alpha_i$.  This equation shows how the real part decreases, on average, with step number. The result has a geometrical interpretation; averaging over diagrams giving zero, the steps with prime number account for \textit{all} of this decrease.
\rm

To motivate the importance of $\zeta (s)$, the sum of terms $cos(\alpha_i logx)$ was shown as a function of $x$.\footnote{Riemann was certainly aware of this representation of primes in 1859, but its properties are now attributed to Landau because of his analysis.}  It can be expressed in the form

\begin{equation}
\sum_{pairs  \quad \alpha_i}^T x^{1/2 \pm i \cdot \alpha_i } = 2 x^{1/2} \sum_{\alpha_i}^T cos(\alpha_i log x)
\end{equation}

 for zeros $1/2 + i \cdot \alpha_i$ whose amplitude is $ \le T$.  Landau proved that the sum is equal to

\begin{equation}
- \frac{T}{2\pi} \Lambda (x) + \mathcal{O} (log T)
\end{equation}

where $\Lambda (x) = log(p)$ if x is $p$ or a power $p^k$ of a prime $p$, and is zero otherwise.  The large difference between $Tlogx/{2\pi}$ when x is a prime or prime power, and $logT$ when it is not, was shown in Figure 3.

Consider the addition of real components of all Argand diagrams giving $\zeta = 0$ for $t < T$;

\begin{equation}
\sum_{\alpha_i < T} \sum_{n=1}^{[\alpha_i/\pi]-1/2}n^{-1/2}cos(\alpha_i log(n))
\end{equation}

The Landau sum shows that only terms with $n=p$ are appreciable.  RH is assumed, and each sum is finite so no conditions are placed on the exchange of summation order.  Figure \ref{primesteps} shows the summed real components of the Argand diagrams for $\alpha_i < 1200$, of which there are 813.  The positions of the primes are marked in red.

\begin{figure}[pt]
\centering
\includegraphics[angle=0,width=.6\textwidth]{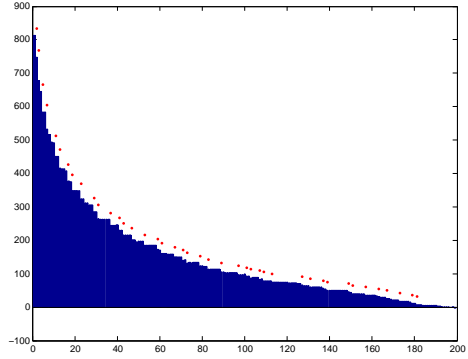}
\caption{Sum of the real parts for Argand diagrams of the 813 zeros of $\zeta$ with $\alpha_i < 1200$.  Primes indicated by red dots.}
\label{primesteps}
\end{figure}

It is also possible to average over a small range of zeros.  Consider the real component $cos(-t log(n))/n^{1/2}$ of steps $n$ for an ensemble of Argand diagrams, each of which gives $\zeta = 0$ for $\alpha_i \in (T, T+ \Delta T)$.  Assume that $\Delta T / T$ is small, so that each diagram has $\simeq T/\pi$ steps. From the classical determination of the number of zeros $\alpha_i<T$, there are about

\begin{equation}
\Delta N = \frac{\Delta T}{2 \pi} log\frac{T}{2 \pi}
\end{equation}
\noindent
zeros on this interval.  Approximating the Landau sum, divided by $x$, by an integral gives the sum of the real components of all of these $\Delta N$ ``zero-Argand'' diagrams.

\begin{equation}
-\frac{\Delta T}{2\pi} \sum_{p=2}^{T/\pi} \frac{logp}{p} \simeq -\frac{\Delta T}{2\pi}\int_2^{T/\pi} \frac{logn}{n} \frac{dn}{logn}
\end{equation}

The sum over primes has been converted to a sum over all integers using the density of primes $logn/n$.  The resulting decrease is $-\Delta N$.  Each real component began at $1$ and ended at $0$, and the Landau sum shows that, on average, the sum of the \textit{prime steps alone} accounts for the total decrease in the real component of $\zeta$ from step 2 to the center of the final spiral.

This process can also be applied to determine the average value of the real part after N steps.  For example, $P(s)=\sum_1^{n_p} n^{-(1/2+\alpha_i)}$ leads to a zero of $\zeta$.  Evaluating the sum in Eqn. 37 to $n_p=\sqrt(T/2 \pi)$ gives $\simeq -\Delta N/2$. This shows that averaging an ensemble of Argand diagrams gives the additional condition that the \textit{real component} of $P(s)$ leading to a zero is $\simeq 1/2$.
\

\it
\item Lehmer's Phenomenon.
\rm

One consequence of the Riemann Hypothesis is based on the product form of $\zeta$ in terms of its zeros.  It states that $Z'/Z$ is monotonic between zeros, which implies that it is not possible for $Z$ to approach zero without changing sign, for example having a local minimum with a positive value.  In geometrical terms, this applies to $\Theta - \phi$, as seen in Figure \ref{phiminthet}. In some of the earliest calculations using digital computers, Lehmer noted places where this forbidden behavior nearly occurred, leading to a possible contradiction with RH\cite[Sec.~8.3.]{edwards}.  One of the occurrences of \textit{Lehmer's Phenomenon} is illustrated in Figure \ref{6704}, near the center of the front edge, which is between Gram points 6707 and 6708.  The lima\c con representation of this point in Figure \ref{limacon6707} demonstrates the ``nonzero close approach'' shows this in terms of the Riemann-Siegel equation:

\begin{figure}[pt]
\centering
\includegraphics[angle=0,width=.6\textwidth]{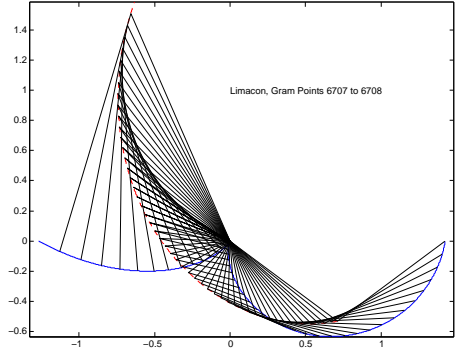}
\caption{Lima\c con for $\sigma = 1/2$, Gram points $6707 \to 6709$.}
\label{limacon6707}
\end{figure}

The dashed red line is the locus of $P(s)$ on this interval.  Radii extend from the origin to points on $P(s)$, from each of which its conjugate radius ends at the solid blue line representing $\zeta (s)$.  Because $\sigma = 1/2$, the lengths of each connected pair of lines are equal.  At this scale, it appears that the $\zeta$ function reaches a cusp at the origin.  The enlarged version in Figure \ref{biglim6707} reveals, however, that there is a very small ``loop'' creating two very closely spaced zeros.\footnote{The points of the curve were computed at much smaller $\delta t$ than the lima\c con lines.}

\begin{figure}[pt]
\centering
\includegraphics[angle=0,width=.6\textwidth]{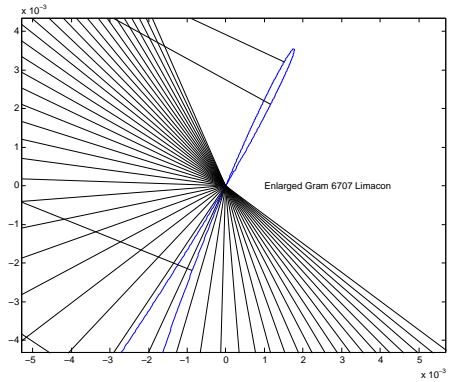}
\caption{Enlarged center portion of Figure \ref{limacon6707}.}
\label{biglim6707}
\end{figure}

This ``close approach'' led Lehmer to speculate that under slightly different conditions, there could have been an extremum which did not lead to a zero, in violation of RH.  Edwards goes on to state that this ``must give pause to even the most convinced believer in the Riemann hypothesis''\cite[Sec.~8.3]{edwards}.
\it
\item Other $\zeta$ Functions
\rm

\begin{figure}[pt]

\centering
\includegraphics[angle=0,width=.6\textwidth]{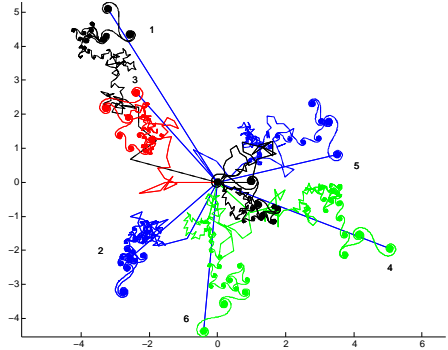}
\caption{The Hurwitz functions $\zeta(s,j/k)$ for the $k=7$ Dirichlet L-functions of $s=.5+10911.9951i$.}
\label{dir7}
\end{figure}

The application of geometrical methods to the Hurwitz $\zeta$ and Dirichlet L-functions shows that they can also be represented by the Argand sum with geometric continuation. These functions display Euler spirals which are not necessarily equal in length to their corresponding initial steps.  The \textit{nth} step has length $1/(n+a)^{-\sigma}$ and angle $-tlog(n+a)$, with angle differences correspondingly defined.   Therefore, the Argand diagrams do not display the detailed symmetry, even for $\sigma = 1/2$, unless $a$ is $0$ or an integer.  There is a final spiral, whose center is the continued value, based on the analysis above.

To illustrate, the set of Dirichlet L-functions for $k=7$ combine the characters in  Table \ref{characters} and the Hurwitz zeta functions shown in Figure \ref{dir7},
where $s= .5+10911.9951i$, a zero of $\zeta$.  Because this is a zero, the function $\zeta(s,r/7)$ for $r=7$ ends the origin.

Further research should reveal applications of the geometric methods displayed here to these functions.

\section{Conclusions}

It has been demonstrated that a comprehensive visual, as opposed to symbolic, understanding of $\zeta (s)$ results from: the geometrical analysis of the symmetry of original steps and their conjugate regions; the identification of spiral centers; the location of the symmetry center $P(s)$; the application of Landau's sum; and the representation of the Riemann-Siegel equation by the lima\c con.  A new description of analytic continuation is given, and primes steps are recognized as the ones which facilitate the transition from $1+i0$ at the end of the first step to the origin for $\zeta$-zeros.  Generalized $\zeta$ functions can be examined by the same techniques, with the restriction that they lack the detailed symmetry principle.

It is the author's belief that the Riemann Hypothesis is correct.  Its proof requires demonstrating that the two values $P(s)$ and $n_p^{1-2\sigma}P(1-s)$ with $s = 1/2 \pm \epsilon + it$ cannot have the same amplitude.

For any interesting mathematical or physical concept, multiple approaches to understanding are beneficial, both to facilitate transmission of ideas and to make new discoveries. As Richard Feynman put it, "Every theoretical physicist who is any good knows six or seven different theoretical explanations for exactly the same physics." \textbf {(The Character of Physical Law.)}

\section{Acknowledgements}
Interest in this subject was kindled many years ago by a physics professor who chose \emph{not} to say "This integral is $\pi^4 /15$, you can look it up in the tables".  Instead, he gave the class a lecture on the zeta function, including its relation to primes via the product form, and "the sport mathematicians call analytic continuation".  The intervening years have led to personal investigations of the popular and technical literature. There have been many illuminating conversations  with colleagues from mathematics, physics or engineering who share an interest in this subject.  At the beginning of the study, Michael R. Stamm of the Air Force Research Laboratory asked "What is that symmetric blob in the middle?".  He has provided the most consistent help and encouragement when difficulties arose.  When shown the geometric derivation for final spiral growth as a function of $\sigma$, David Oro of Los Alamos National Laboratory commented "\textit{That's} why you need analytic continuation.". When asked "What do you mean, \textit{generalized} zeta function?", Robert Whitley of the Mathematics Department at UC Irvine was quite courteous, pointing out the appropriate references. David D. Hardin and Sean Walston, both of the Lawrence Livermore National Laboratory, have read early drafts and made valuable comments.  Encouragement to pursue publication of the Argand diagram approach was provided by Jean-Marc L\'evy-Leblond\cite{leblond}.  Mark Coffey, with an extensive publication record in both physics and analytic number theory at the Colorado School of Mines, has shared many preprints and offered suggestions, including the relevance of the Mellin transform and the importance of a general introductory discussion of $\zeta$. The posting of the first 100,000 zeros to the web by Andrew Odlyzko has been extremely useful (http://dtc.umn.edu).  The assistance of William Buttler of the Physics Division at the Los Alamos National Laboratory for mathematical advice, calculations and figures using MATLAB, and especially for LaTeX help, has been absolutely indispensable.

\end{enumerate}
\end{document}